\documentclass{article}
\usepackage{amsmath,mathrsfs,amssymb,amsthm,amscd}
\usepackage[]{fontenc}
\usepackage{hyperref}
\usepackage{graphics}
\usepackage{a4wide}
\usepackage[dvips]{graphicx,color}

\newcommand{\ezeta}{Z}

\newcommand{\Id}{\mathcal{I}}
\newcommand{\HN}{\mathcal{H}}

\newcounter{EQNR}
\renewcommand{\contentsname}{Contents\\
{\footnotesize\normalfont(A table of contents should normally not
be included)}}

\begin{document}

\title{Zeta functions, heat kernels and spectral asymptotics on
degenerating families of discrete tori\thanks{We are grateful
to the Hausdorff Center in Bonn for the opportunity for us
to meet in the summer of 2007. The first-named and second-named
authors acknowledge support from NSF and PSC-CUNY grants.
The last-mentioned author thanks VR, the G\"oran Gustafsson Foundation,
and G.S. Magnussons fond for support.}}
\author{G. Chinta, J. Jorgenson, and A. Karlsson}
\date{22 December 2008 \\[5mm] Revised: 2 November 2009}

\maketitle

\begin{abstract}\noindent
By a discrete torus we mean the Cayley graph associated to a
finite product of finite cycle groups with generating set given by
choosing a generator for each cyclic factor. In this article we
study the spectral theory of the combinatorial Laplacian for
sequences of discrete tori when the orders of the cyclic factors
tend to infinity at comparable rates. First we show that the
sequence of heat kernels corresponding to the degenerating family
converges, after re-scaling, to the heat kernel on an associated
real torus. We then establish an asymptotic expansion, in the
degeneration parameter, of the determinant of the combinatorial
Laplacian. The zeta-regularized determinant
of the Laplacian of the limiting real torus appears as the
constant term in this expansion. On the other hand, using a classical
theorem by Kirchhoff the determinant of the combinatorial Laplacian of a finite
graph divided by the number of vertices equals the number of
spanning trees, called the complexity, of the graph.  As a result,
we establish a precise connection between the complexity of the
Cayley graphs of finite abelian groups and heights of real tori.
It is also known that spectral determinants on discrete tori can
be expressed using trigonometric functions and that spectral
determinants on real tori can be expressed using modular forms on
general linear groups. Another interpretation of our analysis is
thus to establish a link between limiting values of certain
products of trigonometric functions and modular forms. The heat
kernel analysis which we employ uses a careful study of $I$-Bessel
functions.  Our methods extend to prove the asymptotic behavior of
other spectral invariants through degeneration, such as special
values of spectral zeta functions and Epstein-Hurwitz type zeta
functions.
\end{abstract}

\section{Introduction}\label{1}

\begin{nn}\label{1.1}
The problems we study in the present article begin with the
following very elementary question.  For any $d \geq 1$, let $N =
(n_{1}, \cdots, n_{d})$ denote a $d$-tuple of positive
integers, and consider the product
\begin{equation}\label{ddefinition}
D(N) = \prod\limits_{K \neq 0}\left(2d - 2\cos(2\pi k_{1}/n_{1}) -
\cdots - 2\cos(2\pi k_{d}/n_{d})\right);
\end{equation}
where the product is over all $d$-tuples $K = (k_{1}, \cdots,
k_{d})$ of non-negative integers with $k_{j} < n_{j}$, omitting
the zero vector in the product. The basic question is the
following:  What is the asymptotic behavior of $D(N)$ as
$N\rightarrow\infty$? Setting $V(N) = n_{1} \cdots n_{d}$, then
$(\log D(N))/V(N)$ can be interpreted as an improper Riemann sum, and we
have the limiting formula
\begin{equation}\label{improperintegral}
\frac{1}{V(N)}\log D(N) \rightarrow
\int\limits_{{\mathbf Z}^{d}\backslash {\mathbf
R}^{d}}\log\left(2d-2\cos(2\pi x_{1})- \cdots -
2\cos(2\pi x_{d})\right)dx_{1}\cdots dx_{d} \,\,\,\,\,\,\,\textrm{as
each $n_{j} \rightarrow \infty$.}
\end{equation}
The integral in (\ref{improperintegral}) exists in the sense of improper
integrals from elementary calculus, and the convergence of $D(N)/V(N)$ to the
improper integral in (\ref{improperintegral}) can be verified easily the monotonicity
of $\log x$ for $x > 0$ and calculus.

One can view $D(N)$ as a determinant of a naturally defined
matrix from graph theory. Quite generally, associated to any
finite graph, there is a discrete Laplacian which acts on the finite
dimensional space of complex valued functions whose domain of
definition is the space of vertices of the graph.  With our
normalization of the Laplacian, defined in section \ref{2.4} below,
$D(N)$ is equal to the product of the non-zero eigenvalues of the
Laplacian associated to a graph which we call a discrete torus.

The study of $D(N)$ takes on an entirely new level of significance
beginning with the 1847 paper of Kirchhoff \cite{Kirchhoff} which
further recognizes $D(N)$ as a fundamental invariant from graph
theory. A spanning tree of a graph is a sub-graph which contains
precisely one path between each pair of vertices of the original
graph. Kirchhoff's Theorem states that the number of spanning
trees is equal to $D(N)/V(N)$, the product of the non-zero
eigenvalues of the Laplacian divided by the number of vertices.

Modern mathematics, theoretical computer science, and statistical
physics contain numerous studies which in some way involves the
number of spanning trees of a given graph, determinants of
Laplacians, or other symmetric functions of eigenvalues which form
basic invariants. With all this, we see that the invariant $D(N)$
has considerable significance far beyond the elementary
considerations which allow for its definition.
\end{nn}

\begin{nn}\label{1.2}
\textbf{Summary of the main results.} Our work begins with the
results from \cite{KN} which establishes a theta inversion formula
for the discrete Laplacian acting on the space of vertices of
$n{\mathbf Z} \backslash {\mathbf Z}$. Let $I_{x}(t)$ be the
classical $I$-Bessel function, reviewed in detail in section
\ref{2.1} below. Then using the general concept of theta functions
with inversion formulas constructed from heat kernels (see, for
example, \cite{JLang1}), the authors in \cite{KN} prove that for
any $t > 0$ and integer $x$, the following identity holds for the
theta function associated to the discrete torus $ n{\mathbf Z} \backslash
{\mathbf Z}$:
\begin{equation}\label{thetainversion}
\theta_{n}(t,x) = \frac{1}{n}\sum\limits_{k=0}^{n-1} e^{-(2-2\cos(2\pi
k/n))t + 2\pi i k x/n} = e^{-2t} \sum\limits_{j =
-\infty}^{\infty}I_{x+jn}(2t),
\end{equation}
where $I_x(t)$ denotes the $I$-Bessel function, see section \ref{2.2}
The generalization of (\ref{thetainversion}) to the
$d$-dimensional discrete torus, defined as the product space
$$
DT_{N} = \prod\limits_{j=1}^{d} n_{j}{\mathbf Z} \backslash
{\mathbf Z},
$$
comes from taking a $d$-fold product of the theta functions in
(\ref{thetainversion}).  By computing an
integral transform of the $d$-dimensional theta inversion, we
obtain the following formula, given in Theorem \ref{3.6}. For any
$s \in \mathbf C$ with $\text{\rm Re}(s^{2}) > 0$, we have that
\begin{equation}\label{ehzeta}
\sum\limits_{\Lambda_{j} \neq 0}\log
\left(s^{2}+\Lambda_{j}\right) = V(N)\Id_d(s) + \HN_{N}(s)
\end{equation}
where
$$
\{\Lambda_{j}\}=\{2d - 2\cos(2\pi k_{1}/n_{1}) - \cdots -
2\cos(2\pi k_{d}/n_{d})\} \,\,\,\,\,\textrm{\rm with $k_{j}=0,
\cdots, n_{j}-1$},
$$
as in (\ref{ddefinition}),
$$
\Id_d(s)=\int\limits_{0}^{\infty}
\left(e^{-2dt}e^{-s^{2}t}I_{0}(2t)^{d}-e^{-t}\right)\frac{dt}{t},
$$
and
$$
\HN_{N}(s) = - \int\limits_{0}^{\infty}\left(e^{-s^{2}t}
\left[\theta_{N}(t) -V(N)  e^{-2dt}I_{0}(2t)^{d}- 1\right]+
e^{-t}\right)\frac{dt}{t}.
$$

We now consider a sequence of integral vectors $N(u) = (n_1(u),...,n_{d}(u))$
parameterized by $u \in \mathbf Z$ such that $n_j(u)/u
\rightarrow \alpha_{j}$ as $u\rightarrow \infty$ for each $j$.
We let $A$ be the diagonal matrix with the numbers $\alpha_i$ on
the diagonal
and let $V(A) = \alpha_{1}
\cdots \alpha_{d} \neq 0$.  Through a careful study of the
infinite series of $I$-Bessel functions in (\ref{thetainversion}),
we obtain the following theorem, which is one of the main results
of the present article.

\textbf{Main Theorem.} \emph{Let $\log\det^{\ast}{\mathbf
\Delta}_{DT,N(u)}$ be the log-determinant of the Laplacian of
non-zero eigenvalues on the $d$-dimensional discrete torus
associated to $N(u)$, and $\log\det^{\ast}{\mathbf
\Delta}_{RT,A}$ be the log-determinant of the Laplacian on the
real torus $A{\mathbf Z}^{d} \backslash {\mathbf R}^{d}$.  Then}
\begin{equation}\label{maintheorem}
\log\det\,\!\!^{\ast}{\mathbf \Delta}_{DT,N(u)} =
V(N(u))\mathcal{I}_{d}(0) + \log u^{2} +
\log\det\,\!\!^{\ast}{\mathbf \Delta}_{RT,A} +
o(1)\,\,\,\,\,\,\,\textrm{as $u \rightarrow \infty$}
\end{equation}
\emph{where}
$$
\mathcal{I}_{d}(0)  = \left(\log 2d -
\int\limits_{0}^{\infty}\left(e^{-2dt}(I_{0}(2t)^{d}-1\right)\frac{dt}{t}\right).
$$

There are a number of facets of (\ref{maintheorem}) which are
interesting.  First, the determinant $\det^{\ast}{\mathbf
\Delta}_{DT,N(u)}$ is a legitimate finite product of eigenvalues,
whereas $\det^{\ast}{\mathbf \Delta}_{RT,A}$ is defined through
zeta function regularization.  From this point of view, we have
connected a zeta regularized determinant with a classical
determinant. Moreover, since the first terms in
(\ref{maintheorem}) are universal, this allows for the possibility
of transfering knowledge, for example, between the minimal
spectral determinants for real tori, as studied in
\cite{Chiu,SaSt}, to minimal spectral determinants for discrete
tori, or alternatively, the minimal number of spanning trees. In
this context it is relevant to remark that the main theorem holds
without any changes for general tori, see section \ref{6.8}.
Second, the discrete torus $DT_{N}$ can be viewed as a lattice on
the real torus $A{\mathbf Z}^{d} \backslash {\mathbf R}^{d}$, and
the degeneration which occurs when $u$ tends to infinity amounts
to considering a family of discrete tori which are becoming
uniformly dense in $A{\mathbf Z}^{d} \backslash {\mathbf R}^{d}$.
The main theorem above proves a type of ``re-scaled continuity''
when studying the asymptotic behavior of spectral determinants.
Finally, the classical Kronecker limit formula for Epstein zeta
functions amounts to the evaluation of $\log\det^{\ast}{\mathbf
\Delta}_{RT,A}$ in terms of a generalization of Dedekind's eta
function to a $\textrm{GL}(d, {\mathbf Z})$ modular form.  With
this, we have established a precise connection between the
asymptotics of the number of spanning trees on families of
discrete tori and modular forms.

The verification of the main theorem in the
case $d=1$ can be carried out directly;  the details are presented
in section (\ref{6.0}). When $d=2$, additional explicit
computations are possible, in which case the main theorem becomes
the following result.

\textbf{Main Theorem in the case $d=2$.} \emph{Let $N(u)=(n_1(u),
n_2(u))=(n_1, n_2)$, and assume that $N(u)/u \rightarrow
(\alpha_{1}, \alpha_{2})$.  With the above notation, we have, for
any integer $K > 3$, the asymptotic formula}
\begin{equation}\label{maintheoremd=2}
\log\det\,\!\!^{\ast}{\mathbf \Delta}_{DT,N(u)} =
n_1n_2\frac{4G}{\pi} + \log (n_1n_2) +
\log (|\eta(i\alpha_2/\alpha_1)|^4\alpha_2 /\alpha_1) +
\sum_{k=1}^{K-3}F_k(u)+O(u^{-K})
\end{equation}
\emph{as $u \rightarrow \infty$, where $G$ is the classical
Catalan constant, $\eta$ denotes Dedekind's eta function, and the
functions $F_k(u)$ are explicitly computable and satisfy the
asymptotic bound $F_{k}(u)=O(u^{-k-2})$ as $u \rightarrow \infty$.
}

The Main Theorem (\ref{maintheorem}) comes from studying the
special value $s=0$ in (\ref{ehzeta}).  The analysis we develop in
the proof of (\ref{maintheorem}) extends to prove the asymptotic
behavior of (\ref{ehzeta}) for all $s \in \mathbf C$ with
$\textrm{Re}(s^{2})
> 0$ as well as for special values of the spectral zeta function
$\sum_{j}\Lambda_{j}^{-w}$ for $w \in \mathbf C$ with
$\textrm{Re}(w) > 0$. In subsection \ref{6.5}, we prove the
following theorem.

\textbf{Theorem.} \label{theoremspecialvalues} \emph{With notation
as above, let $\zeta_{A}$ be the spectral zeta function on the
real torus $A{\mathbf Z}^{d} \backslash {\mathbf R}^{d}$.  Then
for any $w\in \mathbf C$ with $\textrm{\rm Re}(w)
> 0$, we have that}
\begin{equation}
\label{sconvergence2} \lim\limits_{u \rightarrow \infty}
u^{-2w}\left(\zeta_{N}(w) -
\frac{V(N)}{\Gamma(w)}\int\limits_{0}^{u^{2}}\left(e^{-2t}(I_{0}(2t))
\right)^{d}t^{w}\frac{dt}{t}\right) = \zeta_{A}(w)
-\frac{V(A)}{(4\pi)^{d/2}(w-d/2)\Gamma(w)}.
\end{equation}
\emph{In particular, for any $w\in \mathbf C$ with $\textrm{\rm
Re}(w)
> d/2$, we have}
\begin{equation}\label{sconvergence1}
\lim\limits_{u \rightarrow \infty}\left(u^{-2w}\zeta_{N(u)}(w)\right)
=\zeta_{A}(w).
\end{equation}

\end{nn}

\begin{nn}\label{1.4}
\textbf{Comparisons with known results.} The case $d=2$ has
been studied since at least the 1960s by physicists,
starting by Kasteleyn's celebrated computation \cite{Kasteleyn}
of the lead term $4G/\pi$.
The next order terms in the asymptotic expansion was later obtained by
Barber, see Duplantier-David \cite{DD} for references and \cite[(3.18)]{DD} for the precise
statements.  Going further, in \cite{DD} the authors establish the
error term $O(u^{-1})$, while our explicit error term expansion
begins a term of order $O(u^{-3})$. The authors in \cite{DD} prove
their results by studying the asymptotic behavior of the spectrum
of the Laplacian, as opposed to the asymptotic behavior of the
heat kernel and its integral transforms, which is the approach
taken in the present article. Subsequent authors have studied
asymptotics of the Laplacian on more general subgraphs of
${\mathbf Z}^2$, see, for example, Burton-Pemantle \cite{BP} and
Kenyon \cite{Kenyon}. For general $d \geq 2$, the articles
\cite{BP, ChangShrock, FelkerLyons, ShrockWu, SokalStarinets} have
discussions which give the lead term asymptotics in
(\ref{maintheorem}).  However, there is no discussion in these
papers or elsewhere regarding the next order term in the
asymptotic expansion in (\ref{maintheorem}).

As stated, in sections \ref{6.5} and \ref{6.6} we apply our
analysis to study the asymptotic behavior of the families of
spectral zeta functions and Epstein-Hurwitz zeta functions through
degeneration.  These problems are considered in \cite{DD} in the
case $d=2$ for general arguments for the Epstein-Hurwitz zeta
function and for the special values $w=1$ and $w=2$ of the
spectral zeta function (see section \ref{6.6} for a clarification
of the notation).  Specifically, equation (3.24) of \cite{DD}
asserts that, ``after a rather long algebra" (quoting \cite{DD}),
one has the asymptotic formula

\begin{equation}
\begin{array}{ll}\label{DDmoment} \displaystyle
\sum_{k\neq0}\frac{1}{\Lambda_{k}^{2}}  & \displaystyle \sim
\left( \frac{\alpha _{1}u\alpha_{2}u}{2\pi}\right)
^{2}\Bigg(\frac{1}{2^{4}\cdot45}(2\pi\alpha
_{1}/\alpha_{2})^{2}+\frac{\zeta(3)}{2(2\pi\alpha_{1}/\alpha_{2})}\\[5mm]
& \displaystyle +\frac{1}{2\pi\alpha_{1}/\alpha_{2}}\sum_{n\geq1}\frac{1}{n^{3}}%
\frac{e^{-2\pi\alpha_{1}/\alpha_{2}n}}{1-e^{-2\pi\alpha_{1}/\alpha_{2}n}}%
+\sum_{n\geq1}\frac{1}{n^{2}}\frac{e^{-2\pi\alpha_{1}/\alpha_{2}n}%
}{(1-e^{-2\pi\alpha_{1}/\alpha_{2}n})^{2}}\Bigg)
\end{array}
\end{equation}
as $u \rightarrow \infty$. For comparison, our result
(\ref{sconvergence1}) specialized to $d=w=2$, yields the
asymptotic formula
\begin{equation}\label{momentd2w2}
\sum_{k\neq0}\frac{1}{\Lambda_{k}^{2}}\sim\frac{u^{4}}{(2\pi)^{4}}
\sum_{(n,m)\neq(0,0)}\frac{1}{((n/\alpha_{1})^{2}+(m/\alpha_{2})^{2})^{2}}
\end{equation}
as $u \rightarrow \infty$.   Going further, the current literature
contains only a handful of other considerations beyond
(\ref{DDmoment}), all in the case when $d=2$, whereas our result
is completely general.

As an aside, note that by equating the lead terms in
(\ref{DDmoment}) and (\ref{momentd2w2}) and setting
$y=\alpha_{1}/\alpha_{2}$, we arrive at the formula
\begin{equation}
\begin{array}{ll}\label{moment2identity} \displaystyle
\sum_{(n,m)\neq(0,0)}\frac{1}{(n^{2}+(my)^{2})^{2}} &
\displaystyle = \left( \frac{2\pi}{y}\right)
^{2}\Bigg(\frac{(2\pi y)^{2}}{2^{4}\cdot45}+\frac{\zeta(3)}{4\pi y}\\[5mm]
& \displaystyle +\frac{1}{2\pi y}\sum_{n\geq1}\frac{1}{n^{3}}%
\frac{e^{-2\pi y n}}{1-e^{-2\pi yn}}%
+\sum_{n\geq1}\frac{1}{n^{2}}\frac{e^{-2\pi yn}%
}{(1-e^{-2\pi yn})^{2}}\Bigg).
\end{array}
\end{equation}

A direct proof of (\ref{moment2identity}) can be obtained by
taking the Fourier expansion of the non-holomorphic Eisenstein
series $E(z,s)$ for $\textrm{\rm SL}(2,{\mathbf Z})$ with $z=iy$
and $s=2$, using the evaluation of the $K$-Bessel function
$K_{s-1/2}$ in terms of exponential functions when $s=2$. We thank
Professor Cormac O'Sullivan for clarifying this point for us, and
we refer the interested reader to his forthcoming article
(\cite{OS}) for a systematic presentation of identities of this
form.

In Riemannian geometry, the determinant of the Laplacian obtained
through zeta function regularization has been used extensively as
a height function on moduli space; see, for example, \cite{Chiu}
and \cite{OPS}.  We view the asymptotic expansion in
(\ref{maintheorem}) as establishing a precise connection to the
well established notion of complexity from graph theory.  Compare
with the discussion in \cite[p. 619]{Sarnak} or \cite[p.
242]{Kenyon}.  Another context where the number of spanning trees
in discrete tori appears is in the study of sandpile models
\cite{CDVV}. Going beyond the setting of discrete and real tori,
we are optimistic that the point of view taken in the present
paper will extend to address problems in the study of other
families of Cayley graphs of discrete, possibly infinite, groups.

\end{nn}

\begin{nn}\label{1.3}
\textbf{Outline of the paper.} In section \ref{2}, we establish
notation and present background material from elsewhere in the
mathematical literature.

In section \ref{3}, we study the theta function
(\ref{thetainversion}) associated to the action of the discrete
Laplacian on functions on $DT_{N} = \prod\limits_{j=1}^{d}
n_{j}{\mathbf Z} \backslash {\mathbf Z}$.  Following the results
from \cite{KN}, we begin with the inversion formula
(\ref{thetainversion}) obtained by expressing the heat kernel on
${\mathbf Z}^{d}$ in terms of $I$-Bessel functions and then
periodize to obtain an expression for the heat kernel on the
discrete torus. The main result in section \ref{3} is Theorem
\ref{3.6}, stated in (\ref{ehzeta}) above, which computes the
Gauss transform (Laplace transform with a quadratic change of
variables) of the heat kernel on $DT_{N}$.  By taking a special
value of the Gauss transform, we obtain an expression for the
determinant of the Laplacian on $DT_{N}$ in terms of integral
transforms of $I$-Bessel functions.

In section \ref{4}, we prove general bounds for $I$-Bessel
functions, building from the fundamental estimates proved in
\cite{Paltsev}. Many of the computations in section \ref{4}
involve re-scaled $I$-Bessel functions.  For instance,
we establish asymptotic behavior
(Proposition \ref{4.7}) and uniform bounds (Lemma \ref{4.6}) in
the parameter $u$ of
\begin{equation}\label{rescalingI}
f(u; x,t) = ue^{-u^{2}t}I_{ux}(u^{2}t)
\end{equation}
In particular, we prove that as $u$ tends to infinity, then $f(u;x, t)$
approaches the heat kernel on ${\mathbf R}$ associated to a
certain scaling of the usual Laplacian.

In section \ref{5} we define and study spectral asymptotics on
degenerating sequences of discrete tori.  Theorem \ref{5.8}, and
its reformulation in Remark \ref{5.9}, proves the main result of
this article as stated in (\ref{maintheorem}), which is the
asymptotic behavior of the determinant of the Laplacian on a
degenerating family of discrete tori. The two key steps are to
first employ the change of variable $t\rightarrow u^2t$ and then,
afterwards, utilize a careful decomposition of the Bessel
integrals involved. The various integrals are then analyzed
individually and matched up with terms in the meromorphic
continuation of the spectral determinant on the continuous side.

Whereas we highlight Theorem
\ref{5.8} as the main result of the present article, there are
many other applications of the analysis from sections \ref{3} and
\ref{4} going beyond the study of the determinant of the Laplacian
which is the focus of attention in section \ref{5}.

In section \ref{6}, we deduce our Main Theorem for $d=2$ and show
that it agrees with the work in \cite{DD}. Moreover, we show how
to obtain an explicit expansion of the error term. This does not
appear in previous works and could be of importance when comparing
heights of discrete tori and corresponing continuous tori.

The main points we address in section \ref{7} are the following:
Examination of the numerical evaluation for the lead term in
various asymptotic expansions for arbitrary dimension;
determination of asymptotic behavior of spectral zeta functions
and Epstein-Hurwitz zeta functions for general dimension;
 and investigation of the notational changes
needed to consider general sequences of degenerating discrete
tori.
\end{nn}

\section{Preliminary material}\label{2}

The purpose of this section is to establish notation and recall
relevant material from elsewhere in the mathematical literature.

\begin{nn}\label{2.1}
\textbf{Basic notation.} We use $t \in \mathbf{R}^{+}$ to denote a
positive real variable and $x \in \mathbf{Z}$ to denote an
integer variable.  For any integer $n \geq 1$, we call the
quotient space $n \mathbf{Z} \backslash \mathbf{Z}$ a discrete
circle, and a finite product of discrete circles will be called a
discrete torus.  The product of $d$ discrete circles formed with
the integers $n_{1}, \cdots, n_{d}$ will be denoted by $DT_{N}$
where $N = (n_{j})_{j=1,...,d}$.

For any function $f : \mathbf{Z} \rightarrow \mathbf{R}$, we
define the Laplacian $\Delta_{\mathbf{Z}}$ by
$$
\Delta_{\mathbf{Z}} f(x) = 2f(x) -
\left(f(x+1)+f(x-1)\right).
$$
The Laplacian on $\mathbf{Z}^{d}$ is the sum of $d$
Laplacians $\Delta_\mathbf{Z}$, one for each coordinate.
The spectrum of the Laplacian
$\Delta_{\mathbf{Z}}$ acting on function on $n \mathbf{Z}
\backslash \mathbf{Z}$ is easily computable (see, for example,
\cite{KN}); hence the spectrum on $DT_{N}$ is simply the set of
sums of eigenvalues for each discrete circle.  Specifically,
let $\{\Lambda_{j}\}$ denote the set of eigenvalues of
$\Delta$ acting on function on $DT_{N}$.  Then, with
our normalization of the Laplacian,
$$
\{\Lambda_{j}\} = \{2d-2\cos(2\pi k_{1}/n_{1})- \cdots - 2\cos(2\pi
k_{d}/n_{d})\}
$$
with $k_{1} = 0, \cdots, n_{1}-1$, $\cdots$, $k_{d} = 0, \cdots,
n_{d}-1$. Let $V(N) = n_{1} \cdots n_{d}$, which can be viewed as
a volume of $DT_{N,d}$.  There are $V(N)$ eigenvalues of the
Laplacian on $DT_{N}$, and $V(N)-1$ of the eigenvalues are
non-zero.

The discrete torus $DT_{N}$ gives rise to a graph by inserting an
edge between two points which have a single component that differs by
one (i.e., nearest neighbor).  This is the Cayley graph of the group
 $\prod\limits_{j=1}^{d} n_{j}{\mathbf Z} \backslash
{\mathbf Z}$ with respect to the generators corresponding to the
standard basis vectors of ${\mathbf Z}^{d}$.
 A spanning tree is a subgraph such that given
any two vertices in $DT_{N}$, there is precisely one path within
the subgraph that connects the two vertices.  The Matrix-Tree
theorem (see e.g. \cite{GR}) asserts that
$$
\#\textrm{\ of spanning trees} = \frac{1}{V(N)}
\prod\limits_{\Lambda_{j}\neq 0}\Lambda_{j}.
$$
We will use the notation
$$
\mbox{$\det^{\ast}$}{\mathbf \Delta}_{DT,N} =
\prod\limits_{\Lambda_{j}\neq 0}\Lambda_{j}
$$
where $\det^{\ast}{\mathbf \Delta}_{DT,N}$ denotes the
determinant of the Laplacian on the discrete torus
$DT_{N}$ omitting the zero eigenvalue.
\end{nn}

\begin{nn}\label{2.2}
\textbf{The $I$-Bessel function.} Classically, the $I$-Bessel
function $I_{x}(t)$ is defined as a solution to the differential
equation
$$
t^{2}\frac{d^{2}w}{dt^{2}} + t\frac{dw}{dt} -(t^{2}+x^{2}) = 0.
$$
For integer values of $x$, it is immediately shown that
$I_{x}=I_{-x}$ and, for positive integer values of $x$, we have
the series representation
\begin{equation}\label{Iseries}
I_{x}(t) =
\sum\limits_{n=0}^{\infty}\frac{(t/2)^{2n+x}}{n!\,\Gamma(n+1+x)}
\end{equation}
as well as the integral representation
\begin{equation}\label{Iintegral}
I_{x}(t) =
\frac{1}{\pi}\int\limits_{0}^{\pi}e^{t\cos(\theta)}\cos(\theta
x)d\theta.
\end{equation}
The mathematical literature contains a vast number of articles and
monographs which study the many fascinating properties and
manifestations of the $I$-Bessel functions, as well as other
Bessel functions.  As demonstrated in the analysis in \cite{KN},
basic to our considerations is the relation
$$
I_{x+1}(t) + I_{x-1}(t) = 2\frac{d}{dt}I_{x}(t),
$$
which easily can be derived from the integral representation and
simple trigonometric identities.
\end{nn}

\begin{nn}\label{2.3}
\textbf{Universal bounds for the $I$-Bessel function.} The
elementary definition of the $I$-Bessel function leads to a number
of precise expressions for $I_{x}(t)$, two of which are stated in
section \ref{2.2}.  Unfortunately, the explicit identities do not
easily lead to viable estimates for $I_{x}(t)$.  Beginning with
the differential equation which characterizes $I_{x}(t)$, the
author in \cite{Paltsev} was able to derive very precise upper
and lower bounds for the $I$-Bessel function, which we now state.
Let
$$
g_{x}(t) = \sqrt{(x^{2}+t^{2})} + x\log \left(\frac{t}{x +
\sqrt{(x^{2}+t^{2})}}\right).
$$
Then for all $t > 0$ and $x \geq 2$ we have that
\begin{equation}\label{paltsevbound}
e^{-1/(2\sqrt{(x^{2}+t^{2})})} \leq \sqrt{(2\pi
)}\cdot(x^{2}+t^{2})^{1/4}I_{x}(t)e^{-g_{x}(t)} \leq
e^{1/(2\sqrt{(x^{2}+t^{2})})}.
\end{equation}
\end{nn}
These bounds will play an important role in this article.  Indeed,
in section 4, we will study the function $g_{x}(t) - t$ in order
to obtain elementary bounds for $\sqrt{t}e^{-t}I_{x}(t)$ which
will be vital when establishing our main result.

\begin{nn}\label{2.4}
\textbf{Heat kernels and zeta functions.} In \cite{KN}, the
authors established a theta inversion formula obtained from
studying the heat kernel associated to the Laplacian on
$n\mathbf{Z}\backslash \mathbf{Z}$.  Let $K_{\mathbf{Z}}(t,x)$ be
a function on $\mathbf{Z}$ such that
$$
\Delta_{\mathbf{Z}}K_{\mathbf{Z}}(t,x) +
\frac{\partial }{\partial t}K_{\mathbf{Z}}(t,x)=0
$$
and
$$
\lim\limits_{t \rightarrow 0} K_{\mathbf{Z}}(t,x) = 1
\,\,\,\,\,\textrm{if $x = 0$, and equal to zero if $x \neq 0$.}
$$
General results from \cite{DM} and \cite{Dodziuk} prove the
existence and uniqueness of $K_{\mathbf{Z}}(t,x)$ and, as shown,
for example in \cite{KN}, we have that
$$
K_{\mathbf{Z}}(t,x) = e^{-2t}I_{x}(2t).
$$
The existence and uniqueness theorems for heat kernels on graphs
apply to heat kernels on $n\mathbf{Z}\backslash \mathbf{Z}$.
Bounds for the $I$-Bessel function from section \ref{2.3} imply
that the series
\begin{equation}\label{heatkerneldc}
\sum\limits_{k=-\infty}^{\infty}e^{-2t}I_{kn + x}(2t)
\end{equation}
converges, and hence is equal to the heat kernel on
$n\mathbf{Z}\backslash \mathbf{Z}$, which we denote by
$K_{n\mathbf{Z}\backslash \mathbf{Z}}(t,x).$
Denote
the eigenvalues on $n\mathbf{Z}\backslash \mathbf{Z}$ and
corresponding eigenfunctions by
$\{\lambda_{n,j}\}$ and $\{\phi_{n,j}\}$ respectively.  This
leads to the identity
\begin{equation}\label{thetainversiondc}
\sum\limits_{j=0}^{n-1}e^{-t\lambda_{n,j}}\phi_{n,j}(x)
\phi_{n,j}(0) =\sum\limits_{k=-\infty}^{\infty}e^{-2t}I_{kn +
x}(2t).
\end{equation}
In the case $n=1$, this yields a classical identity which expresses an
infinite sum of $I$-Bessel functions as an exponential function.
The heat kernel on $DT_{N}$ and corresponding theta inversion
formula is obtained by taking the $d$-fold product of
(\ref{thetainversiondc}).
\end{nn}

\begin{nn}\label{2.5}
\textbf{Spectral analysis on real tori.} Given a positive definite
$d \times d$ matrix $A$, let $RT_{A}$ denote the real torus
$A\mathbf{Z}^{d} \backslash \mathbf{R}^{d}$.  On $\mathbf{R}$
we use the variable $y$ to denote the standard global coordinate,
from which we have the Laplacian
$$
\Delta_{\mathbf{R}}f(y) = -\frac{d^{2}}{dy^{2}}f(y)
$$
and corresponding heat kernel $K_{\mathbf{R}}(t,y)$ which is
uniquely characterized by the conditions
$$
\Delta_{\mathbf{R}}K_{\mathbf{R}}(t,y) + \frac{\partial}{\partial
t} K_{\mathbf{R}}(t,y)=0
$$
and
$$
\int\limits_{\mathbf{R}}K_{\mathbf{R}}(t,x-y)f(y)dy \rightarrow
f(x) \,\,\,\,\,\textrm{as $t \rightarrow 0$}
$$
for any smooth, real-valued function $f$ on $\mathbf{R}$ with
compact support.  The set of eigenvalues of the Laplacian on
$RT_{A}$ is given by $\{(2\pi)^{2} \,^{t}mA^{\ast}m\}$ for $m\in
\mathbf{Z}^{d},$ where $A^{\ast}$ is the dual lattice to $A$.  The
eigenfunctions of the Laplacian are expressible as exponential
functions.  The theta function
$$
\Theta_{A}(t)= \sum\limits_{m \in
\mathbf{Z}^{d}}e^{-(2\pi)^{2}\cdot \,^{t}mA^{\ast}m\cdot t},
$$
with $t > 0$, is the trace of the heat kernel on $RT_{A}$.  The
asymptotic behavior of the theta function is well-known, namely
that
$$
\Theta_{A}(t) = V(A)(4\pi t)^{-d/2} + O(e^{-c/t})
\,\,\,\,\,\,\,\textrm{for some $c > 0$ as $t \rightarrow 0$,}
$$
and
$$
\Theta_{A}(t) = 1 + O(e^{-ct}) \,\,\,\,\,\,\,\textrm{for some $c
> 0$ as $t \rightarrow \infty$.}
$$
\end{nn}

\begin{nn}\label{2.6}
\textbf{Regularized determinants on real tori.}
 In very general circumstances, the spectral zeta function is defined
 as the Mellin transform $\mathbf M$ of the theta function formed with the non-zero
eigenvalues. Specifically, for $s \in \mathbf C$ with
$\textrm{Re}(s)
> d/2$, we define the spectral zeta function
\begin{equation}\label{zetaint}
\zeta_{A}(s)={\mathbf M}\Theta_{A}(s) =
\frac{1}{\Gamma(s)}\int\limits_{0}^{\infty}
\left(\Theta_{A}(t)-1\right)t^{s}\frac{dt}{t}.
\end{equation}
The integral in (\ref{zetaint}) converges for $\textrm{Re}(s) >
d/2$.  Let us write
\begin{equation}\label{zetaint2}
\begin{array}{ll}
\displaystyle \zeta_{A}(s) & \displaystyle =
\frac{1}{\Gamma(s)}\int\limits_{0}^{1}\left(\Theta_{A}(t)-V(A)
(4\pi t)^{-d/2}\right)t^{s}\frac{dt}{t} \\[5mm] & \displaystyle +
\frac{1}{\Gamma(s)}\int\limits_{0}^{1}\left(V(A)(4\pi
t)^{-d/2}-1\right)t^{s}\frac{dt}{t}
\\[5mm]&\displaystyle +
\frac{1}{\Gamma(s)}\int\limits_{1}^{\infty}\left(\Theta_{A}(t)-1\right)t^{s}\frac{dt}{t}.
\end{array}
\end{equation}

\noindent Going further, we can carry out the integral in the
second term to get
\begin{equation}\label{zetaint3}
\begin{array}{ll}\displaystyle
\displaystyle \zeta_{A}(s) &\displaystyle =
\frac{1}{\Gamma(s)}\int\limits_{0}^{1}\left(\Theta_{A}(t)-V(A)(4\pi
t)^{-d/2}\right)t^{s}\frac{dt}{t} \\[5mm] &\displaystyle +
(4\pi)^{-d/2}\frac{V(A)}{(s-d/2)\Gamma(s)} - \frac{1}{\Gamma(s+1)}
 \\[5mm]&\displaystyle +
\frac{1}{\Gamma(s)}\int\limits_{1}^{\infty}\left(\Theta_{A}(t)
-1\right)t^{s}\frac{dt}{t}
\end{array}
\end{equation}

\noindent By virtue of the asymptotic behavior of $\Theta_{A}(t)$
as $t \rightarrow 0$ and $t\rightarrow \infty$, the expression
(\ref{zetaint3}) provides a meromorphic continuation of
$\zeta_{A}(s)$ to $s \in \mathbf{C}$. Therefore, using
(\ref{zetaint3}), we can study the behavior of $\zeta_{A}(s)$ near
$s=0$. The integrals are holomorphic near $s=0$, and $1/\Gamma(s)
= s + O(s^{2})$ near $s=0$, so then we have that $\zeta_{A}(0) =
-1$, which is a point that will be used later.  In particular, we
have that
\begin{equation}\label{zetaderivzero}
\begin{array}{ll}\displaystyle
\displaystyle \zeta_{A}'(0) &\displaystyle =
\int\limits_{0}^{1}\left(\Theta_{A}(t)-V(A)(4\pi
t)^{-d/2}\right)\frac{dt}{t} + \Gamma'(1) \\[5mm] &\displaystyle
-\frac{2}{d}V(A) (4\pi)^{-d/2}  +
\int\limits_{1}^{\infty}\left(\Theta_{A}(t)-1\right)\frac{dt}{t},
\end{array}
\end{equation}
\end{nn}
which will play an important role in the proof of our Main Result.

\begin{nn}\label{2.7}
\textbf{Kronecker's limit formula.} The special value
$\zeta_{A}'(0)$ has a unique place in classical analytic number
theory, and the evaluation of $\zeta_{A}'(0)$ in
terms of modular forms is generally referred to as Kronecker's
limit formula, which we now describe.  The content of this section
has its origins in \cite{Epstein}, and the discussion we present
here comes directly from \cite{DI}.

Let $Q$ denote a $d \times d$ positive definite matrix, and $u\in
{\mathbf R}^{d}$ any vector.  If we write $Q =
\left(q_{i,j}\right)$, then the quadratic form associated to $Q$
is defined by
$$
Q(u) = \sum\limits_{i,j}q_{i,j}u_{i}u_{j}.
$$
The Epstein zeta function associated to $Q$ is defined for
$\textrm{\rm Re}(s) > d/2$ by the convergent series
$$
\ezeta(s,Q) = \sum\limits_{m \in \mathbf{Z}^{d}
\setminus\{0\}}Q(m)^{-s}.
$$
From the discussion in sections \ref{2.4} and \ref{2.5}, the
Epstein zeta function is, up to a multiplicative factor of
$(2\pi)^{-2s}$, the spectral zeta function associated to the real
torus $RT_{A}$ with $A = Q^{\ast}$, where $Q^{\ast}$ is the dual
lattice to $Q$. It can be shown that $\ezeta(s,Q)$ admits a
meromorphic continuation and the functional equation
$$
\pi^{-s}\Gamma(s)\ezeta(s, Q^{-1}) = (\det
Q)^{1/2}\pi^{s-d/2}\Gamma(d/2-s)\ezeta(d/2-s,Q).
$$
The Iwasawa decomposition of $Q$ asserts that $Q$ can be uniquely
expressed as
$$
Q = \left(\begin{matrix}1 & 0 \\-^{t}x & I_{d-1} \end{matrix}\right)
\left(\begin{matrix}y^{-1} & 0 \\0 & Y \end{matrix}\right)
\left(\begin{matrix}1 & -x \\0 & I_{d-1} \end{matrix}\right)
$$
with $y \in \mathbf{R}_{+}$, $x \in \mathbf{R}^{d-1}$ and $Y$ a
$(d-1)\times (d-1)$ positive definite matrix; here, $I_{d-1}$
denotes the $(d-1)\times (d-1)$ identity matrix and we think of $x
\in \mathbf{R}^{d-1}$ and $y \in \mathbf{R}$ as row vectors.  For
$m \in \mathbf{R}^{d-1}$, one defines
$$
Q\{m\} = m \cdot x + i \sqrt{yY(m)},
$$
which is a complex number with $\textrm{\rm Im}(Q\{m\}) > 0$
unless $m=0$.  With all this, the generalization of the classical
Kronecker limit formula for $\ezeta(s,Q)$ is the identity
\begin{equation}\label{kronecker}
\ezeta'(0,Q) = -2\pi\sqrt{y}\ezeta(1/2,Y) - \log \left|(2\pi)^{2}y
\prod\limits_{m \in \mathbf{Z}^{n-1}\setminus \{0\}} \left(1 -
e^{2\pi i Q\{m\}}\right)^{4}\right|.
\end{equation}
In \cite{DI} the authors study $\exp(-\ezeta'(0,Q))$ as a
generalization of the classical Dedekind eta function and prove an
analogue of the Chowla-Selberg formula which relates
$\exp(-\ezeta'(0,Q))$ to special values of higher order Gamma
functions and special values of certain $L$ functions.  We refer
the interested reader to \cite{DI} for additional results
regarding the fascinating number theoretic and automorphic aspects of
(\ref{kronecker}).
\end{nn}

\begin{nn}\label{2.8}
\textbf{Mellin transform and inversion.} Let
$f_x(t)=e^{-t/2}I_{x}(t/2).$ Known transform identities include
the formula
\begin{equation}
  \label{eq:3a}
\tilde f_x(s)=\int_0^\infty e^{-t/2}I_{x}(t/2)t^s \frac{dt}{t}
=\frac{\Gamma(s+x)\Gamma(1/2-s)}{\sqrt{\pi}\Gamma(x+1-s)} .
\end{equation}
The identity (\ref{eq:3a}) is valid for $-x<\textrm{Re}(s)<1/2.$
The inverse Mellin formula gives
\begin{equation}
  \label{eq:4a}
  f_x(t)=\frac{1}{2\pi i} \int_{(\sigma)} \tilde f_x(s) t^{-s} ds,
\text{\ \ for \ \ }-x<\sigma <1/2.
\end{equation}

The Mellin transform of a product is computed via the convolution
of two transforms, yielding the computations
\begin{align*}
  \widetilde{f_xf_y}(s) &=\int_0^\infty f_x(t)f_y(t)t^s
\frac{dt}{t}\\
&= \int_0^\infty f_x(t)\left[ \frac{1}{2\pi i} \int_{(\sigma)}
\tilde
  f_y(z) t^{-z} dz \right]t^s \frac{dt}{t}, \,\,\,\,\,\textrm{for $-y<\sigma<1/2$}\\
&= \frac{1}{2\pi i}  \int_{(\sigma)} \tilde f_y(z)\left[
\int_0^\infty f_x(t) t^{s-z}\frac{dt}{t}\right] dz ,
\,\,\,\,\,\textrm{for $-x<\textrm{Re}(s)
-\sigma<1/2$} \\
& = \frac{1}{2\pi i}  \int_{(\sigma)} \tilde f_x(s-z)\tilde f_y(z)
dz.
\end{align*}

Therefore, we conclude that
\begin{equation}
\begin{array}{ll}
  \label{eq:5a}
\displaystyle \widetilde{f_xf_y}(s) & \displaystyle =
\int_0^\infty e^{-t}I_{x}(t/2)I_{y}(t/2)t^s \frac{dt}{t} \\[5mm]&\displaystyle =
\frac{1}{2\pi i}\int_{(\sigma)} \frac{\Gamma(s-z+x)\Gamma(1/2+
  z-s)\Gamma(z+y)\Gamma(1/2-z)}
{\sqrt{\pi}\Gamma(x+1+z-s)\sqrt{\pi}\Gamma(y+1-z)} dz
\\[5mm]&\displaystyle =\frac{1}{2\pi i}\int_{(\sigma)} \frac{\Gamma(s-z+x)\Gamma(1/2+
  z-s)\Gamma(z+y)\Gamma(1/2-z)}
{\pi\Gamma(x+1+z-s)\Gamma(y+1-z)} dz .
\end{array}
\end{equation}
Using the bounds $I_x(t)= O(t^x)$ as $t\rightarrow 0$ and
$e^{-t}I_x(t) = O(t^{-1/2})$ as $t \rightarrow \infty$, the Mellin
transform $\widetilde{f_xf_y}(s)$ is defined for
$-x-y<\textrm{Re}(s) <1$, so (\ref{eq:5a}) holds for
$0<\sigma<\textrm{Re}(s)<1/2.$   If $x\neq 0$ or $y \neq 0$, we
can take $s=0$ and write
\begin{equation}
  \label{eq:6a}
  \int^\infty_0 e^{-t}[I_{x}(t/2) I_{y}(t/2)]\frac{dt}{t} =\frac{1}{2\pi i}\int_{(1/4)}
\frac{\Gamma(x-z)\Gamma(1/2+z)\Gamma(z+y)\Gamma(1/2-z)}
{\pi\Gamma(x+1+z)\Gamma(y+1-z)} dz.
\end{equation}
In section \ref{6.2} below, we will extend (\ref{eq:6a}) to the
case $x=y=0$.
\end{nn}

\begin{nn}\label{2.9}
\textbf{Miscellaneous results.} For any $\varepsilon > 0$ and $z
\in \mathbf C$, Stirling's formula for the classical Gamma
function is the asymptotic relation
$$
\log \Gamma (z) = (z-1/2)\log z - z + (1/2)\log (2\pi) + O(1/z)
\,\,\,\,\,\textrm{as $z \rightarrow \infty$ provided $\vert \arg
(z) \vert < \pi -\varepsilon$.}
$$
In particular, one has for fixed $a$ and $b$, the asymptotic
relation
\begin{equation}\label{gammaratio}
\frac{\Gamma(z+a)}{\Gamma(z+b)} = z^{a-b} +
O\left(z^{a-b-1}\right)\,\,\,\,\,\textrm{as $z \rightarrow \infty$
provided $\vert \arg (z) \vert < \pi -\varepsilon$.}
\end{equation}
Further terms in the Stirling's formula can be computing (see, for
example, \cite{JLang0}), which then would imply further terms in
the asymptotic expansion (\ref{gammaratio}).

Throughout our work, we will use the elementary identity
\begin{equation}
  \label{eq:log}
\log (w) = \int\limits_{0}^{\infty}(e^{-t}-e^{-wt})\frac{dt}{t}
\,\,\,\,\,\textrm{for all $w \in \mathbf{C}$ with $\textrm{\rm
Re}(w) > 0$.}
\end{equation}

To prove this relation, one simply observes that both sides of the
proposed identity vanish when $w=1$ and have first derivative
equal to $1/w$.
\end{nn}

\section{Zeta functions and determinants for discrete
tori}\label{3}

In this section we study the Gauss transform of the trace of the
heat kernel associated to a discrete torus.  We begin with Lemma
\ref{3.1} which recalls the theta inversion formula associated to
a general discrete torus. From the theta inversion formula, we
define the Gauss transform (Lemma \ref{3.2}).  Using the group
periodization representation of the heat kernel, we decompose the
theta inversion formula into two summands:  The identity term and
the set of non-identity terms.  We study the Gauss transform of
these two summands separately, ultimately arriving at an identity
(Theorem \ref{3.6}) which expresses the determinant of the
Laplacian on the discrete torus in terms of $I$-Bessel functions.

\begin{nn}\label{3.1}
\textbf{Lemma.} {\it Let $\theta_{N}(t)$ be the theta
function associated to the $d$-dimensional discrete torus
$N\mathbf{Z}^{d}\backslash \mathbf{Z}^{d} $,
defined by}
$$
\theta_{N}(t) = \sum\limits_{\Lambda_{j}}e^{-\Lambda_{j}t}
$$
\emph{where}
$$
\{\Lambda_{j}\} = \{2d- 2\cos(2\pi m_{1}/n_{1})-\cdots - 2\cos(2\pi
m_{d}/n_{d}): 0\leq m_i< n_i, \ \text{for each $i=1,2,\ldots,d $} \}.
$$
\emph{Then for all $t > 0$, we have the identity}
$$
\theta_{N}(t)= V(N)\sum\limits_{K \in {\mathbf
Z}^{d}}\prod\limits_{1\leq j\leq d}e^{-2t}I_{n_{j}\cdot
k_{j}}(2t),
$$
where $K$ runs over $d$-tuples of integers $(k_1, \ldots, k_d).$
\end{nn}

\proof We refer to \cite{KN} for the case $d=1$.  From there, one gets the
general case using that the heat kernel on a product space is equal
to the product of the heat kernels, since the Laplacian on
the product space is defined to be the sum of the Laplacians from
each factor space.  \hfill $\Box$

\begin{nn}\label{3.2}
\textbf{Lemma.} \emph{For all $s \in {\mathbf C}$ with $\textrm{\rm
Re}(s^{2})> 0$, we have}
$$
\sum\limits_{\Lambda_{j}\neq 0}\frac{2s}{s^{2}+\Lambda_{j}} =
V(N) 2s\int\limits_{0}^{\infty}
e^{-s^{2}t}e^{-2dt}(I_{0}(2t))^{d}dt + 2s\int\limits_{0}^{\infty}
e^{-s^{2}t}\Big[\theta_{N}(t) -V(N)  e^{-2dt}I_{0}(2t)^{d}-
1\Big]dt.$$
\end{nn}
\proof By the definition of $\theta_{N}$, we can write
$$
\sum\limits_{\Lambda_{j}\neq 0}e^{-\Lambda_{j}t} = V(N)
\left(e^{-2t}I_{0}(2t)\right)^{d} +  [\theta_{N}(t) -V(N)
e^{-2dt}I_{0}(2t)^{d}- 1].
$$
Now simply multiply both sides of this identity by $2s
e^{-s^{2}t}$ and integrate with respect to $t$ on $(0,\infty)$.
Asymptotic behavior of the integrands as $t \rightarrow 0$ and
$t\rightarrow \infty$ easily imply that the resulting integrals
are convergent for $s \in \mathbf{C}$ provided $\textrm{\rm
Re}(s^{2})> 0$.
 \hfill $\Box$

\begin{nn}\label{3.3}
\textbf{Lemma} \emph{The function}
$$
f(s) = \sum\limits_{\Lambda_{j} \neq 0}\log
\left(s^{2}+\Lambda_{j}\right)
$$
\emph{is uniquely characterized by the differential equation}
$$
\partial_{s}f(s) = \sum\limits_{\Lambda_{j}\neq 0}\frac{2s}{s^{2}+\Lambda_{j}}
$$
\emph{and the asymptotic relation}
$$
f(s) = (V(N)-1)\cdot \log s^{2} + o(1) \,\,\,\,\,\emph{as $s
\rightarrow \infty$.}
$$
\end{nn}
\proof The differential equation characterizes $f(s)$ up to an
additive constant, which is uniquely determined by the stated
asymptotic behavior.   \hfill $\Box$

\begin{nn}\label{3.4}
\textbf{Proposition}
 \emph{The function}
$$
\Id_d(s)=
-\int\limits_{0}^{\infty}\left(e^{-s^{2}t}e^{-2dt}I_{0}(2t)^{d}
-e^{-t}\right)\frac{dt}{t}
$$
\emph{is uniquely characterized by the differential equation}
$$
\partial_{s}\Id_d(s) =  2s\int\limits_{0}^{\infty}
e^{-s^{2}t}e^{-2dt}I_{0}(2t)^{d}dt
$$
\emph{and the asymptotic relation}
$$
\Id_d(s) =  \log s^{2} + o(1) \,\,\,\,\, \emph{as $s \rightarrow
+\infty$.}
$$
\end{nn}
\proof For this, we write
$$
\Id_d(s)=
-\int\limits_{0}^{\infty}e^{-s^{2}t}e^{-2dt}\left((I_{0}(2t))^{d}-1\right)
\frac{dt}{t} -
\int\limits_{0}^{\infty}\left(e^{-s^{2}t}e^{-2dt}-e^{-t}\right)\frac{dt}{t}.
$$
By equation (\ref{eq:log}),
$$
\int\limits_{0}^{\infty}\left(e^{-s^{2}t}e^{-2dt}-e^{-t}\right)\frac{dt}{t}
= -\log (s^{2}+2d),
$$
so then we have
$$
\Id_d(s)=
-\int\limits_{0}^{\infty}e^{-s^{2}t}e^{-2dt}\left(I_{0}(2t)^{d}-1\right)
\frac{dt}{t}+ \log (s^{2}+2d).
$$
From this last expression, the asymptotic behavior as $s
\rightarrow \infty$ is immediate. \hfill $\Box$

\begin{nn}\label{3.5}
\textbf{Proposition} \emph{The function}
$$
\HN_{N}(s)= -
\int\limits_{0}^{\infty}\left(e^{-s^{2}t}\left[\theta_{N}(t)
-V(N)  e^{-2dt}I_{0}(2t)^{d}- 1\right]+ e^{-t}\right)\frac{dt}{t}
$$
\emph{is uniquely characterized by the differential equation}
$$
\partial_{s}\HN_{N}(s) =   2s\int\limits_{0}^{\infty}
e^{-s^{2}t}\Big[\theta_{N}(t) -V(N)  e^{-2dt}I_{0}(2t)^{d}-
1\Big]dt
$$
\emph{and the asymptotic relation}
$$
\HN_{N}(s) = - \log s^{2} + o(1) \,\,\,\,\, \emph{as $s
\rightarrow +\infty$.}
$$
\end{nn}
\proof For this, we write
\begin{equation}\label{hintegral}
\begin{array}{ll}
\displaystyle  \HN_{N}(s) &\displaystyle  = -
\int\limits_{0}^{\infty} e^{-s^{2}t}\left(\theta_{N}(t) -V(N)
e^{-2dt}I_{0}(2t)^{d}\right) \frac{dt}{t} +
\int\limits_{0}^{\infty}(e^{-s^{2}t}-e^{-t})\frac{dt}{t}
\\[7mm]& \displaystyle
= - \int\limits_{0}^{\infty} e^{-s^{2}t}\left(\theta_{N}(t)
-V(N) e^{-2dt}I_{0}(2t)^{d}\right) \frac{dt}{t} - \log(s^{2}).
\end{array}
\end{equation}
Clearly, the integral in (\ref{hintegral}) approaches zero as $s$
approaches infinity, which completes the proof of the stated
asymptotic relation. \hfill $\Box$

\begin{nn}\label{3.6}
\textbf{Theorem.} \emph{For any $s\in {\mathbf C}$ with $\textrm{\rm
Re}(s^{2}) > 0$, we have the relation}
$$
\sum\limits_{\Lambda_{j} \neq 0}\log
\left(s^{2}+\Lambda_{j}\right) = V(N)\Id_d(s) + \HN_{N}(s).
$$
\emph{Letting $s\to 0$, we have the identity}
$$
\log \left(\prod\limits_{\Lambda_{j}\neq 0}\Lambda_{j}\right) =
V(N)\Id_d(0) + \HN_{N}(0)
$$
\emph{where}
$$
\Id_d(0) =
-\int\limits_{0}^{\infty}\left(e^{-2dt}I_{0}(2t)^{d}-e^{-t}\right)\frac{dt}{t}
$$
\emph{and}
$$
\HN_{N}(0)  = -
\int\limits_{0}^{\infty}\left(\theta_{N}(t) -
V(N)e^{-2dt}I_{0}(2t)^{d} - 1 + e^{-t}\right)\frac{dt}{t}.
$$
\end{nn}

\proof We begin with the relation from Lemma \ref{3.2}.  By
substituting from the differential equations from Lemma \ref{3.3},
Proposition \ref{3.4} and Proposition \ref{3.5} and then
integrating, we get that
\begin{equation}\label{gausstrans}
\sum\limits_{\Lambda_{j} \neq 0}\log
\left(s^{2}+\Lambda_{j}\right) = V(N)\Id_d(s) + \HN_{N}(s) +
C
\end{equation}
for some constant $C.$  We now use the asymptotic behavior relations
as $s \rightarrow \infty$ from Lemma \ref{3.3}, Proposition
\ref{3.4} and Proposition \ref{3.5} to show that $C=0$.  From the
series expansion (\ref{Iseries}), we have that
$$
e^{-2dt}I_{0}(2t)^{d}-e^{-t} = O(t) \,\,\,\,\,\,\,\textrm{as $t
\rightarrow 0$.}
$$
Lemma \ref{4.1} gives that
$$
e^{-2dt}I_{0}(2t)^{d} =
O\left(t^{-d}\right)\,\,\,\,\,\,\,\textrm{as $t \rightarrow
\infty$,}
$$
so then the integrand in the definition of $\Id_d(0)$ is
$L^{1}(0,\infty)$ with respect to $dt/t$.  Concerning the
integrand in the definition of $\HN_N(0)$ we have that
$$
\theta_{N}(t) - V(N)e^{-2dt}(I_{0}(2t))^{d} = O(t)
\,\,\,\,\,\,\,\textrm{as $t \rightarrow 0$}
$$
and
$$
e^{-t} - 1 = O(t) \,\,\,\,\,\,\,\textrm{as $t \rightarrow 0$}
$$
so then the integrand in the definition of $\HN_N(0)$ is in
$L^{1}(0,1)$ with respect to $dt/t$. Furthermore, we have
$$
\theta_{N}(t) - 1 = O(e^{-ct}) \,\,\,\,\,\,\,\textrm{as $t
\rightarrow \infty$ for some $c > 0$}
$$
and
$$
V(N)e^{-2dt}(I_{0}(2t))^{d} - e^{-t} = O(t^{-d})
\,\,\,\,\,\,\,\textrm{as $t \rightarrow \infty$,}
$$
so then the integrand in the definition of $\HN_N(0)$ is in
$L^{1}(1,\infty)$ with respect to $dt/t$. With all this, we have
that all functions in (\ref{gausstrans}) are continuous and
well-defined for $s \in {\mathbf R}_{\ge 0}$, so we simply need to
evaluate at $s=0$ to complete the proof.

 \hfill $\Box$

\section{Bounds and asymptotic formulas for Bessel
functions}\label{4}

In this section we prove bounds for individual $I$-Bessel
functions. For technical reasons, it is necessary to separately
establish bounds for $I_{0}$ and for $I_{x}$ for $x > 0$.  As will
be seen in the next section, it is necessary to consider the
function
$$
u\cdot e^{-u^{2}t}I_{n(u)}(u^{2}t)
$$
for $u \geq 1$ and where $n(u)/u\rightarrow x$. The asymptotic
behavior as $u\rightarrow \infty$ is given in Proposition
\ref{4.7}, uniform bounds are established in Lemma \ref{4.1} when
$x=0$, and in Lemma \ref{4.4} and Lemma \ref{4.6} when $x > 0$.
Although we restrict our attention to integers $x \geq 0$, we
recall that $I_{x} = I_{-x}$, so the results we prove apply for
all $x \in \textbf{N}$.

\begin{nn}\label{4.1}
\textbf{Lemma.} \emph{For any $\varepsilon < \pi /2$ and $t > 0$,
we have the bounds}
$$
0 \leq e^{-t}I_{0}(t) \leq C\cdot t^{-1/2}\,\,\,\,\,
\textrm{\emph{where} $C =\frac{1}{\sqrt{(2-\varepsilon^{2}/6)\pi}}
+ \frac{\pi - \varepsilon}{\pi}\cdot
\frac{1}{\sqrt{(1-\varepsilon^{2}/12)\varepsilon^{2}e}}$.}
$$
\end{nn}
\proof

The positivity of $I_{0}(t)$ follows immediately from the series
expansion stated in section \ref{2.2}.  From the integral
representation, choose any $\varepsilon \in (0,\pi/2)$ and write
$$
e^{-t}I_{0}(t)=
\frac{1}{\pi}\int\limits_{0}^{\varepsilon}e^{-t(1-\cos(u))}du +
\frac{1}{\pi}\int\limits_{\varepsilon}^{\pi}e^{-t(1-\cos(u))}du.
$$
There exists $c(\varepsilon) > 0$ such that $1-\cos(u) \geq
c(\varepsilon)u^{2}$ for $u \in [0,\varepsilon]$.  With this, for
the first integral we have
$$
\frac{1}{\pi}\int\limits_{0}^{\varepsilon}e^{-t(1-\cos(u))}du \leq
\frac{1}{\pi}\int\limits_{0}^{\varepsilon}e^{-c(\varepsilon)tu^{2}}du
\leq \frac{1}{\sqrt{4c(\varepsilon)\pi t}}.
$$
For the second integral, we trivially have
$$
\frac{1}{\pi}\int\limits_{\varepsilon}^{\pi}e^{-t(1-\cos(u))}du
\leq \frac{\pi - \varepsilon}{\pi}e^{-t(1-\cos(\varepsilon))} \leq
\frac{\pi - \varepsilon}{\pi}e^{-c(\varepsilon)t\varepsilon^{2}}.
$$
Combining, we have
$$
0 \leq e^{-t}I_{0}(t) \leq \frac{1}{\sqrt{4c(\varepsilon)\pi t}} +
\frac{\pi - \varepsilon}{\pi}e^{-c(\varepsilon)t\varepsilon^{2}}.
$$
Since $\varepsilon < \pi/2$ and $u \leq \varepsilon$, we have the
bound
$$
1-\cos(u) \geq u^{2}/2 - u^{4}/24 \geq
(1/2-\varepsilon^{2}/24)u^{2},
$$
so we may take $c(\varepsilon) = (1/2-\varepsilon^{2}/24)$.
Therefore,
$$
0 \leq e^{-t}I_{0}(t) \leq \frac{1}{\sqrt{(2-\varepsilon^{2}/6)\pi
t}} + \frac{\pi -
\varepsilon}{\pi}e^{-(1/2-\varepsilon^{2}/24)\varepsilon^{2}\cdot
t}.
$$
Using elementary calculus, one shows that
$$
t^{1/2}e^{-at} \leq 1/\sqrt{2ae},
$$
so then
\begin{equation}\label{Ibound}
0 \leq e^{-t}I_{0}(t) \leq C\cdot t^{-1/2}\,\,\,\,\, \textrm{where
$C =\frac{1}{\sqrt{(2-\varepsilon^{2}/6)\pi}} + \frac{\pi -
\varepsilon}{\pi}\cdot
\frac{1}{\sqrt{(1-\varepsilon^{2}/12)\varepsilon^{2}e}}$.}
\end{equation}
 \hfill $\Box$

\begin{nn}\label{4.2}
\textbf{Remark.} Directly from Lemma \ref{4.1}, we have that
$$
u e^{-u^{2}t}I_{0}(u^{2}t) \leq u \cdot C\cdot (u^{2}t)^{-1/2} =
C\cdot t^{-1/2},
$$
so we indeed have established the uniform upper bound as claimed.
Also, for our purposes, it is not necessary to optimize
(\ref{Ibound}) through a judicious choice of $\varepsilon$.
Numerically, one can show that by taking $\varepsilon = 1.5$,
which is allowed since we only required that $\varepsilon <
\pi/2$, we have that $C=0.676991\dots$.  One point we will use
later (see section \ref{6.8a}) is that $C < 1/\sqrt{2} < 1$.  We
will use the numerical verification of this bound and omit the
theoretical proof from our analysis, noting that the estimate for
$C$ indeed can be proved from (\ref{Ibound}).
\end{nn}

\begin{nn}\label{4.3}
\textbf{Lemma.} \emph{For fixed $x \geq 1$ and $t>0$, consider the
function}
$$
h_{x}(t) = \sqrt{(x^{2}+t^{2})}-t+x\log
\left(\frac{t}{x+\sqrt{(x^{2}+t^{2})}}\right).
$$
\emph{Then we have the bound}
$$
\exp(h_{x}(t)) \le \left(\frac{t}{t+x}\right)^{x/2}.
$$
\emph{By continuity, the inequality also holds when $t=0$.}
\end{nn}
\proof To begin, observe that
\begin{equation}\label{expon2}
h_{x}(t) = \frac{x^{2}}{t+\sqrt{(x^{2}+t^{2})}}+x\log
\left(\frac{t}{x+\sqrt{(x^{2}+t^{2})}}\right),
\end{equation}
which comes from the definition of $h_{x}(t)$ and by writing
$$
\sqrt{(x^{2}+t^{2})}-t = \frac{x^{2}}{t+\sqrt{(x^{2}+t^{2})}}.
$$
Therefore, we have that
$$
\begin{array}{ll}\displaystyle
h_{x}(t) &\displaystyle =
\frac{x^{2}}{t+\sqrt{(x^{2}+t^{2})}}+x\log
\left(\frac{t}{x+\sqrt{(x^{2}+t^{2})}}\right)
\\[7mm] & \displaystyle =
\frac{x^{2}}{t+\sqrt{(x^{2}+t^{2})}} + x \log
\left(\frac{1}{t}\left(\sqrt{x^{2}+t^{2}} - x \right)\right)
\\[7mm] & \displaystyle =
\frac{x^{2}}{t+\sqrt{(x^{2}+t^{2})}} + x \log
\left(\sqrt{\left(1+\frac{x^{2}}{t^{2}}\right)} -
\frac{x}{t}\right)
\\[7mm] & \displaystyle = x \left(
\frac{x/t}{1+\sqrt{((x/t)^{2}+1)}} + \log
\left(\sqrt{\left(1+\frac{x^{2}}{t^{2}}\right)} -
\frac{x}{t}\right)\right).
\end{array}
$$
We now employ the change of variables
$$
u =\log \left(\sqrt{\left(1+\frac{x^{2}}{t^{2}}\right)} +
\frac{x}{t}\right)
$$
which is equivalent to the relation $\sinh (u) = x/t$. Using the
elementary identities
$$
1 + (\sinh u)^{2} =(\cosh u)^{2},\,\,\,\,\,\sinh u = 2 \sinh
(u/2)\cosh(u/2)
$$
and
$$
1+\cosh u = 2 (\cosh (u/2))^{2},\,\,\,\,\,\cosh u - \sinh u =
e^{-u}
$$
to arrive at the expression
$$
h_{x}(t) = x \left(\tanh(u/2)-u\right).
$$
Trivially, since $u\geq 0$ we have that $\tanh(u/2) \leq u/2$, so
then, for $x>0$,
$$
h_{x}(t) \displaystyle  \leq  -x\cdot u /2  = -(x/2)\log
\left(\sqrt{\left(1+\frac{x^{2}}{t^{2}}\right)} +
\frac{x}{t}\right)  \leq  -(x/2)\log \left(1 + \frac{x}{t}\right).
$$
With all this, we have that
$$
\exp(h_{x}(t))  \leq \left(1 + \frac{x}{t}\right)^{-x/2}
 = \left(\frac{t}{x+t}\right)^{x/2},
$$
which completes the proof of the lemma. \hfill $\Box$

\begin{nn}\label{4.4}
\textbf{Corollary:} \emph{For any $t > 0$ and integer $x \geq 0$,
we have}
$$
\sqrt{t}\cdot e^{-t} I_{x}(t) \leq
\left(\frac{t}{t+x}\right)^{x/2} =
\left(1+\frac{x}{t}\right)^{-x/2}.
$$
\end{nn}

\proof We begin by considering $x \geq 2$, so then the analysis
from Lemma \ref{4.3} applies.  Indeed, we use the trivial
estimates
$$
\frac{1}{\sqrt{2\pi}}e^{1/(2\sqrt{(x^{2}+t^{2})})} \leq
\frac{e^{1/4}}{\sqrt{2\pi}} \leq 1
$$
and
$$
\frac{\sqrt{t}}{(x^{2}+t^{2})^{1/4}} \leq 1.
$$
Therefore, using the notation of Lemma \ref{4.3}, the bound
(\ref{paltsevbound}) becomes
$$
\sqrt{t}\cdot e^{-t}I_{x}(t) \leq \exp(h_{x}(t)) \leq
\left(\frac{t}{x + t}\right)^{x/2},
$$
which proves the claim, again provided that $x \geq 2$. If $x=0$,
the claim follows from Lemma \ref{4.1} as well as Remark \ref{4.2}
which shows that the constant $C$ in Lemma \ref{4.1} satisfies $C
\leq 1$.  It remains to consider the case when $x=1$.

The series representation (\ref{Iseries}) of $I_{x}(t)$ gives
$$
0 \leq I_{1}(t) \leq (t/2)
\sum\limits_{n=0}^{\infty}\frac{(t/2)^{2n}}{n!\Gamma(n+2)} \leq
(t/2)\left(\sum\limits_{n=0}^{\infty}\frac{(t/2)^{n}}{n!}\right)
\left(\sum\limits_{n=0}^{\infty}\frac{(t/2)^{n}}{\Gamma(n+2)}\right)
\leq (t/2)e^{t}.
$$
For $t \leq 1$, $t/2 \leq 1/\sqrt{(t+1)}$, so then
$$
\sqrt{t}\cdot e^{-t}I_{1}(t) \leq \sqrt{t}\cdot (t/2) \leq
\left(\frac{t}{t+1}\right)^{1/2}\,\,\,\,\,\textrm{for $t \leq 1$.}
$$
From (\ref{Iintegral}), we have that $I_{1}(t) \leq I_{0}(t)$
since $\cos(x\theta) \leq 1$. With this, Lemma \ref{4.1} implies
\begin{equation}\label{I1bound1}
e^{-t}I_{1}(t) \leq e^{-t}I_{0}(t) \leq C/\sqrt{t}
\end{equation}
In Remark \ref{4.2}, it was argued that $C < 1/\sqrt{2}$.  With
this, we have for $t \geq 1$ the inequalities
$$
C/\sqrt{t} \leq 1/\sqrt{2t} \leq 1/\sqrt{(t+1)},
$$
so then
\begin{equation}\label{I1bound2}
e^{-t}I_{1}(t) \leq 1/\sqrt{(t+1)}\,\,\,\,\,\textrm{for $t \geq
1$.}
\end{equation}
Combining (\ref{I1bound1}) and (\ref{I1bound2}), we have the
claimed assertion for $x=1$, which completes the proof for all
integers $x \geq 0$.
 \hfill $\Box$

 \begin{nn}\label{4.5}
\textbf{Remark.} To be precise, the bound in Corollary \ref{4.4}
in the case $x=1$ is not needed in this article.  We included the
statement and proof for the sake of completeness.
\end{nn}

\begin{nn}\label{4.6}
\textbf{Lemma.} \emph{Fix $t \geq 0$ and non-negative integers $x$
and $n_{0}$. Then for all $n \geq n_{0}$, we have the uniform
bound}
$$
0 \leq \sqrt{(n^{2}t)}\cdot e^{-n^{2}t}I_{nx}(n^{2}t) \leq
\left(\frac{n_{0}t}{x+n_{0}t}\right)^{n_{0}x/2} = \left(1 +
\frac{x}{n_{0}t}\right)^{-n_{0}x/2} .
$$
\end{nn}
\proof The positive lower bound is obvious from the series
definition of the $I$-Bessel function, so we focus on proving the
upper bound.  From Corollary \ref{4.4}, we have
$$
\sqrt{(n^{2}t)}\cdot e^{-n^{2}t}I_{nx}(n^{2}t) \leq
 \left(\frac{n^{2}t}{nx + n^{2}t}\right)^{nx/2}.
$$
Elementary algebra yields
$$
\left(\frac{n^{2}t}{nx+n^{2}t}\right)^{-1} = 1 + \frac{x}{nt} = 1
+ \frac{x^{2}/(2t)}{(nx/2)},
$$
so then we have
\begin{equation}\label{Iscalebound}
\sqrt{(n^{2}t)}\cdot e^{-n^{2}t}I_{nx}(n^{2}t) \leq
\left(\frac{n^{2}t}{nx + n^{2}t}\right)^{nx/2} = \left(1 +
\frac{x^{2}/(2t)}{(nx/2)}\right)^{-nx/2}.
\end{equation}
For any constant $c > 0$, consider the function
$$
g(y) = \left(1 + \frac{c}{y}\right)^{y}.
$$
We claim that $g(y)$ is monotone increasing in $y$.  Indeed, using
logarithmic differentiation, we have that
$$
\begin{array}{ll}
\displaystyle \frac{g'(y)}{g(y)} & \displaystyle = \log
\left(1+\frac{c}{y}\right) + y \cdot \frac{1}{1+c/y}\cdot
\frac{-c}{y^{2}} = \log \left(1+\frac{c}{y}\right) - \frac{c}{c+y}
= \log \left(\frac{c+y}{y}\right) - \frac{c}{c+y}
\\[5mm] &\displaystyle = -\log \left(\frac{y}{c+y}\right)- \frac{c}{c+y}
= -\log \left(1-\frac{c}{c+y}\right)- \frac{c}{c+y}.
\end{array}
$$
Clearly, if $0 \leq u < 1$, the function $-\log (1-u)-u$ is
positive.  Since $g(y) > 0$, we conclude that $g'(y) > 0$, which
shows that $g(y)$ is monotone increasing in $y$.  By taking
$c=x^{2}/(2t)$ and $y=nx/2$, we have, for $n \geq n_{0}$, the
inequality
$$
\left(1 + \frac{x^{2}/(2t)}{(n_{0}x/2)}\right)^{n_{0}x/2} \leq
\left(1 + \frac{x^{2}/(2t)}{(nx/2)}\right)^{nx/2}
$$
which gives
\begin{equation}\label{ratebound}
\left(1 + \frac{x^{2}/(2t)}{(nx/2)}\right)^{-nx/2} \leq \left(1 +
\frac{x^{2}/(2t)}{(n_{0}x/2)}\right)^{-n_{0}x/2} = \left(1 +
\frac{x}{n_{0}t}\right)^{-n_{0}x/2} =  \left(
\frac{n_{0}t}{n_{0}t+x}\right)^{n_{0}x/2}.
\end{equation}
By substituting (\ref{ratebound}) in (\ref{Iscalebound}), the
lemma is proved. \hfill $\Box$

\begin{nn}\label{4.7}
\textbf{Proposition.} \emph{Let $n(u)$ be a sequence of positive
integers parameterized by $u \in \mathbf{Z}^{+}$ such that
$n(u)/u\rightarrow \alpha >0$ as $u \rightarrow \infty$.  Then for
any $t > 0$ and non-negative integer $k \geq 0$, we have}
$$
\lim\limits_{u\rightarrow \infty} n(u)\cdot
e^{-2u^{2}t}I_{n(u)k}(2u^{2}t)= \frac{\alpha}{\sqrt{4\pi t}}\cdot
e^{-(\alpha k)^{2}/(4t)}.
$$
\end{nn}
\proof  We work with the integral expression for the $I$-Bessel
function given in section \ref{2.2}, namely that for any integer
$x$, we have
$$
I_{x}(t) =
\frac{1}{\pi}\int\limits_{0}^{\pi}e^{t\cos(\theta)}\cos(\theta
x)d\theta.
$$
To begin, assume $k> 0$. If we let $y=n(u)k\theta$, then we can write
$$
\begin{array}{ll}
\displaystyle n(u)e^{-2u^{2}t}I_{n(u)k}(2u^{2}t) & \displaystyle  =
\frac{n(u)}{\pi}\int\limits_{0}^{\pi}e^{-2u^{2}t(1-\cos(\theta))}\cos(\theta
n(u)k)d\theta
\\[5mm] &\displaystyle =
\frac{1}{k\pi}\int\limits_{0}^{n(u)k\pi}e^{-2u^{2}t(1-\cos(y/(n(u)k)))}\cos(y)dy.
\end{array}
$$
For all $v \in [0,\pi]$, one can easily show that
$$
\frac 12 - \frac{\pi^2}{24}\leq \frac 12 - \frac
{v^2}{24}\leq\frac{1-\cos(v)}{v^2}.
$$
Let $c = 1/2 - \pi^{2}/24$, which, numerically, can be shown to
satisfy $c > 0$.  Setting $v=y/(n(u)x)$, we get for $y\in [0,
ux\pi]$ the uniform bound
$$
\frac{1-\cos(y/(n(u)k))}{1/u^{2}} \geq
c\left(\frac{yu}{n(u)k}\right)^{2}.
$$
In addition, observe that
$$
\lim\limits_{u\rightarrow \infty}\frac{1-\cos(y/(n(u)k))}{1/u^{2}} =
\frac{1}{2} \left(\frac{y}{\alpha k}\right)^{2}.
$$
Choose $u_{0}$ so that for $u > u_{0}$, we have  $\alpha/2 <
n(u)/u < 2\alpha$. Using elementary bounds, we have for any $u >
u_{0}$ the inequalities
$$
\left| n(u)e^{-2u^{2}t}I_{n(u)k}(2u^{2}t) \right| \leq
\frac{1}{k\pi}\int\limits_{0}^{n(u)k\pi}e^{-cy^{2}t/(2\alpha
k)^{2}}dy  \leq
\frac{1}{k\pi}\int\limits_{0}^{\infty}e^{-cy^{2}t/(2\alpha
k)^{2}}dy = \frac{\alpha }{\sqrt{\pi c t}}.
$$
Therefore, by the Lebesgue Dominated Convergence Theorem, we have
$$
\lim\limits_{u\rightarrow \infty} n(u)\cdot
e^{-2u^{2}t}I_{n(u)k}(2u^{2}t)=
\frac{1}{k\pi}\int\limits_{0}^{\infty}e^{-y^{2}t/(\alpha k)^{2}}\cos(y)dy
= \frac{\alpha}{\sqrt{4\pi t}}\cdot e^{-(\alpha k)^{2}/(4t)}.
$$
If $x=0$, the proof follows a similar pattern using instead the
substitution $y=\theta u$; for the sake of
brevity, we omit the details.   \hfill $\Box$

\begin{nn}\label{4.8}
\textbf{Remark.}
Proposition \ref{4.7} was proved by Athreya \cite[Theorem 2]{A} using a
certain local central limit theorem. For convenience we
included a quick self-contained proof.
\end{nn}

\section{Asymptotic behavior of spectral determinants}\label{5}

In this section we use the bounds and asymptotic relations from
the previous section to prove our main theorem, namely the
asymptotic behavior of the spectral determinant associated to a
sequence of degenerating discrete tori.  We begin by proving that
the associated family of traces of heat kernels converges pointwise
through degeneration (Proposition \ref{5.3}) and then prove
uniform bounds for heat traces for long time (Proposition
\ref{5.5}) and small time (Proposition \ref{5.6} and Lemma
\ref{5.7}).  After these results, we will analyze the expression
for $\HN_{N}(0)$ from Lemma \ref{3.5}, namely
\begin{equation}\label{hterm2}
\HN_{N}(0)  = -
\int\limits_{0}^{\infty}\left(\theta_{N}(t) -
V(N)e^{-2dt}I_{0}(2t)^{d} - 1 + e^{-t}\right)\frac{dt}{t}.
\end{equation}
Ultimately we compare the limiting value of (\ref{hterm2}) through
degeneration with the expression for the spectral determinant on
the real torus, which is stated in (\ref{zetaderivzero}) in
section \ref{2.6}. These computations are given in the proof of
Theorem \ref{5.8}. which is the main result of this section.

\begin{nn}\label{5.1}
\textbf{Degenerating sequences of discrete tori.} We consider
$d$-tuples of integers $N(u)$ parametrized by a positive integer $u$
in such a way that
\begin{equation}\label{degen}
\frac{1}{u}N(u) = \frac{1}{u}(n_{1},\cdots, n_{d}) \rightarrow
(\alpha_{1},\cdots, \alpha_{d}) \,\,\,\,\,\,\textrm{as $u
\rightarrow \infty$.}
\end{equation}
Let $A$ be the diagonal matrix with the $\alpha_j$'s on the
diagonal. Recall the notation $V(N) = n_{1}\cdots n_{d}$ and
$V(A)=\alpha_{1}\cdots \alpha_{d}$. From (\ref{degen}),
$V(N(u))/u^{d} \rightarrow V(A)$ when $u \rightarrow \infty$. For
the sake of brevity, and, for this section, we assume that $V(A) >
0$, meaning the limiting real torus $A\mathbf Z^{d}\backslash
\mathbf R^{d}$ has dimension $d$. Using the elementary change of
variables $t \mapsto u^{2}t$, we write (\ref{hterm2}) as
\begin{equation}\label{hterm3}
\HN_{N(u)}(0)  = -
\int\limits_{0}^{\infty}\left(\theta_{N(u)}(u^{2}t) -
V(N(u))\left(e^{-2u^{2}t}I_{0}(2u^{2}t)\right)^{d} - 1 +
e^{-u^{2}t}\right)\frac{dt}{t},
\end{equation}
which is the form that we will study. The results of this section
are designed to determine the asymptotic behavior of
(\ref{hterm3}) as $u \rightarrow \infty$.
\end{nn}

\begin{nn}\label{5.2}
\textbf{Proposition.} \emph{For each fixed $t > 0$, we have the
pointwise limit}
$$
\theta_{N(u)}(u^{2}t) \rightarrow \Theta_A(t) \,\,\,\,\,\textrm{as
$u \rightarrow \infty$.}
$$
\emph{In words, the rescaled theta functions on the discrete tori
approach the theta function on the limiting real torus $A\mathbf
Z^{d}\backslash \mathbf R^{d}.$}
\end{nn}
\proof We begin by writing the theta function using its expansion
involving $I$-Bessel functions, namely
\begin{equation}\label{thetaseries}
\theta_{N(u)}(u^{2}t) = \sum\limits_{k_{1},\cdots, k_{d} =
-\infty}^{\infty}\prod\limits_{j=1}^{d}n_{j}(u)e^{-2u^{2}t}
I_{n_{j}(u)k_{j}}(2u^{2}t).
\end{equation}
From Proposition \ref{4.7}, we have, for any $t > 0$, the
pointwise limit
$$
\prod\limits_{j=1}^{d}n_{j}(u)e^{-2u^{2}t}I_{n_{j}(u)k_{j}}(2u^{2}t)
\rightarrow V(A) \prod\limits_{j=1}^{d}\frac{1}{\sqrt{(4\pi
t)}}e^{-(\alpha_{j}k_{j})^{2}/(4t)}\,\,\,\,\,\,\,\textrm{as $u
\rightarrow \infty$.}
$$
Since
$$
\Theta_A(t) = V(A)\sum\limits_{k_{1},\cdots, k_{d} =
-\infty}^{\infty}\prod\limits_{j=1}^{d}\frac{1}{\sqrt{(4\pi
t)}}e^{-(\alpha_{j}k_{j})^{2}/(4t)},
$$
the result will follow if we can interchange the limit in $u$ with
the infinite sum in (\ref{thetaseries}).

Choose $u_{0}$ such that for each $j = 1, \cdots, d$,
$\alpha_{j}/2 < n_{j}(u)/u < 2 \alpha_{j}$ for all $u > u_{0}$.
We can re-write (\ref{thetaseries}) into a sum of $d+1$ subseries,
each determined by the number of $k_{j}$'s which are equal to
zero. Fix some $u_{0} \geq 0$ sufficiently large. Then from Lemma
\ref{4.1} and Lemma \ref{4.6}, we have that the subseries
consisting of terms with exactly
 $r$ of the $k_j$'s equal to zero is bounded from above for any $u > u_{0}$
 by
\begin{equation}\label{thetaseriesbound}
2^{r}C^{r}(2t)^{-r/2} V(A) \sum\limits_{k_{1},\cdots, k_{r} =
1}^{\infty}\prod\limits_{j=1}^{r} \left(1 +
\frac{\alpha_{j}k_{j}}{4u_{0}t}\right)^{-u_{0}\alpha_{j}k_{j}/2}
\leq 2^{r}C^{r}(2t)^{-r/2} V(A) \sum\limits_{k_{1},\cdots, k_{r} =
1}^{\infty}\prod\limits_{j=1}^{r} r_{j}^{k_{j}}
\end{equation}
where
\begin{equation}\label{ratio}
r_{j} = \left(1 +
\frac{\alpha_{j}}{4u_{0}t}\right)^{-u_{0}\alpha_{j}/2} < 1.
\end{equation}
Obviously,
$$
\sum\limits_{k_{1},\cdots, k_{r} =
1}^{\infty}\prod\limits_{j=1}^{r} r_{j}^{k_{j}} =
\prod\limits_{j=1}^{r} \left(\frac{r_{j}}{1-r_{j}}\right).
$$
Therefore, the series in (\ref{thetaseries}) is uniformly
convergent for fixed $t$, so we can interchange the limit in $u$
and the summation in (\ref{thetaseries}), which completes the
proof. \hfill $\Box$

\begin{nn}\label{5.3}
\textbf{Lemma.} \emph{Let}
$$
\theta_{\textrm{abs}}(t) =  2 \sum\limits_{j=1}^{\infty}e^{- c
j^{2} t}
$$
\emph{with $c = 4\pi^{2}(1-\pi^{2}/24)^{2}$.  Assume $u_{0}$ is
such that $\alpha_{k}/2 < n_{k}/u < 2 \alpha_{k}$ for all $k =
1,\cdots, d$ and $u > u_{0}$.  Then for any $t
> 0$ and $u > u_{0}$, we have the bound }
$$
\theta_{N(u)}(u^{2}t) \leq
\prod\limits_{k=1}^{d}\left(1+e^{-4u_{0}^{2}t} +
\theta_{\textrm{abs}}(\alpha_{k}t/2)\right).
$$
\end{nn}
\proof Recall from \cite{KN} that
$$
\theta_{n_{k}(u)}(t) = 1 + e^{-4t}+ 2
\sum\limits_{j=1}^{n_{k}(u)/2-1}e^{-4(\sin (\pi j/n_{k}(u)))^{2}t}
$$
if $n_{k}(u)$ is even, and
$$
\theta_{n_{k}(u)}(t) = 1  + 2
\sum\limits_{j=1}^{(n_{k}(u)-1)/2}e^{-4(\sin (\pi
j/n_{k}(u)))^{2}t}
$$
if $n_{k}(u)$ is odd. In either case, we will use the elementary
bound
$$
\sin x \geq x - x^{3}/6 \,\,\,\,\,\textrm{ for $x \in [0,\pi/2]$}.
$$
To prove this inequality, one considers the Taylor series
expansion with error at $x=0$ for $\sin (x)$, together with the
observation that the error term is positive.  Taking $x = \pi
j/n_{k}(u)$, we then have that
$$
n_{k}(u) \sin(\pi j/n_{k}(u)) \geq \pi j(1-(\pi
j)^{2}/(6n_{k}(u)^{2})).
$$
We may assume that $j\leq n_{k}(u)/2$, so then $j/n_{k}(u) \leq
1/2$, hence
$$
n_{k}(u) \sin(\pi j/n_{k}(u)) \geq \pi j(1-\pi^{2}/24).
$$
It is important to use the trivial bound $1 -\pi^{2}/24
> 0$. With all this, we arrive at the inequality
$$
\theta_{n_{j}(u)}(u^{2}t) \leq 1 + e^{-4u^{2}t} + 2
\sum\limits_{j=1}^{n/2}e^{- c \alpha_{k}j^{2} t/2} \leq 1 +
e^{-4u_{0}^{2}t} + 2 \sum\limits_{j=1}^{\infty}e^{- c \alpha_{k}
j^{2} t/2}
$$
with $c = 4\pi^{2}(1-\pi^{2}/24)^{2} > 0$.  Therefore, we have
shown that for any $u > u_{0}$, we have the bound
$$
\theta_{n}(u^{2}t) \leq 1 + e^{-4u_{0}^{2}t} +
\theta_{\textrm{abs}}(t).
$$
Since
$$
\theta_{N(u)}(u^{2}t) =
\prod\limits_{j=1}^{d}\theta_{n_{k}(u)}(u^{2}t),
$$
the proof of the lemma is complete.  \hfill $\Box$

\begin{nn}\label{5.4}
\textbf{Proposition.} \emph{With notation as above, we have that }
$$
\begin{array}{ll} \displaystyle \int\limits_{1}^{\infty} &
\displaystyle \left(\theta_{N(u)}(u^{2}t) -
V(N(u))(e^{-2u^{2}t}I_{0}(2u^{2}t))^{d} - 1+
e^{-u^{2}t}\right)\frac{dt}{t}  \\[7mm] & \displaystyle =
\int\limits_{1}^{\infty}\left(\Theta_A(t)-1\right)\frac{dt}{t}
-\frac{2}{d}V(A)(4\pi)^{-d/2} + o(1) \,\,\,\,\,\textrm{as $u
\rightarrow \infty$.}
\end{array}
$$
\end{nn}
\proof  Write
\begin{equation}
\begin{array}{ll} \label{longterm}\displaystyle \int\limits_{1}^{\infty}& \displaystyle
 \left(\theta_{N(u)}(u^{2}t) -
V(N(u))(e^{-2u^{2}t}I_{0}(2u^{2}t))^{d} - 1+
e^{-u^{2}t}\right)\frac{dt}{t} \\ [7mm] & \displaystyle =
\int\limits_{1}^{\infty} \left(\theta_{N(u)}(u^{2}t) -
1\right)\frac{dt}{t}  -  \int\limits_{1}^{\infty} \displaystyle
V(N(u))(e^{-2u^{2}t}I_{0}(2u^{2}t))^{d} \frac{dt}{t} +
\int\limits_{1}^{\infty} e^{-u^{2}t}\frac{dt}{t},
\end{array}
\end{equation}
and we consider the three integrals separately.  Trivially, we
have
$$
\int\limits_{1}^{\infty}e^{-u^{2}t} \frac{dt}{t} \rightarrow 0
\,\,\,\,\,\textrm{as $u \rightarrow \infty$.}
$$
Next, we claim that
\begin{equation}\label{identasymp}
V(N)\int\limits_{1}^{\infty}\left(e^{-2u^{2}t}I_{0}(2u^{2}t)\right)^{d}
\frac{dt}{t} \rightarrow V(A)\int\limits_{1}^{\infty}(4\pi
t)^{-d/2}\frac{dt}{t} = \frac{2}{d}V(A)(4\pi)^{-d/2}
\,\,\,\,\,\textrm{as $u \rightarrow \infty$.}
\end{equation}
Indeed, the pointwise convergence of the integrands in
(\ref{identasymp}) is proved in Proposition \ref{4.7}, and the
uniform bound from Lemma \ref{4.1}, which is integrable on
$(1,\infty)$ with respect to the measure $dt/t$, allows one to
apply the Lebesgue Dominated Convergence Theorem to prove the
claim.  It remains to show that
$$
\int\limits_{1}^{\infty} \left(\theta_{N(u)}(u^{2}t) -1
\right)\frac{dt}{t} \rightarrow
\int\limits_{1}^{\infty}\left(\Theta_A(t)-1\right)\frac{dt}{t}
\,\,\,\,\,\textrm{as $u \rightarrow \infty$.}
$$
The pointwise convergence of the integrands is proved in
Proposition \ref{5.2}, and Lemma \ref{5.3} establishes a uniform,
integrable upper bound so that, again, we may apply the Lebesgue
Dominated Convergence Theorem, thus completing the proof.
 \hfill $\Box$

\begin{nn}\label{5.5}
\textbf{Proposition.} \emph{With notation as above, we have that }
$$
\int\limits_{0}^{1} \left(\theta_{N(u)}(u^{2}t)
-V(N)(e^{-2u^{2}t}I_{0}(2u^{2}t))^{d} \right)\frac{dt}{t}
\rightarrow \int\limits_{0}^{1}\left(\Theta_A(t)-V(A)(4\pi
t)^{-d/2}\right)\frac{dt}{t}
$$
\emph{as $u \rightarrow \infty$.}
\end{nn}

\proof For fixed $t < 1$, we have the pointwise convergence
$$
\theta_{N(u)}(u^{2}t) -V(N)(e^{-2u^{2}t}I_{0}(2u^{2}t))^{d}
\rightarrow \Theta_A(t)-V(A)(4\pi t)^{-d/2}
\,\,\,\,\,\textrm{as $u \rightarrow \infty$}
$$
directly from Proposition \ref{4.7} and Proposition \ref{5.2}.  It
remains to show that we have uniform, integrable bounds for the
integrand so then the proposition will follow from the Lebesgue
Dominated Convergence Theorem.  For this, we write
\begin{equation}\label{smalltime}
\theta_{N(u)}(u^{2}t) -V(N(u))(e^{-2u^{2}t}I_{0}(2u^{2}t))^{d}
= V(A) \sum\limits_{K \neq 0}\prod\limits_{k_{1},\cdots,
k_{d}}ue^{-2u^{2}t}I_{n_{j}(u) k_{j}}(2u^{2}t).
\end{equation}
The bounds (\ref{thetaseriesbound}) and (\ref{ratio}) apply. Since
each $n_{j}(u)$ tends to infinity, one can choose $u_{0}$ so that
for each $j = 1, \cdots, d$ and $u > u_{0}$, we have the
inequality $n_{j}(u)k_{j}
> d +2$, which implies that the upper bound in (\ref{ratio})
is integrable on $(0,1)$ with respect to $dt/t$.
 \hfill $\Box$

\begin{nn}\label{5.6}
\textbf{Lemma.} \emph{For $u \in \mathbf R$, we have the asymptotic
formula}
$$
\int\limits_{0}^{1}\left(e^{-u^{2}t}-1\right)\frac{dt}{t}=\Gamma'(1)
- \log (u^{2}) +o(1) \,\,\,\,\,\textrm{as $u \rightarrow \infty$.}
$$
\end{nn}

\proof Choose $\varepsilon \in (0,1)$.  By employing integration
by parts, we can write
$$
\begin{array}{ll}\displaystyle
\int\limits_{\varepsilon}^{1}\left(e^{-u^{2}t}-1\right)\frac{dt}{t}
& \displaystyle =
\int\limits_{\varepsilon}^{1}\left(e^{-u^{2}t}-1\right)d\log (t) =
\log (t)\cdot\left(e^{-u^{2}t}-1\right)\Big|_{t=\varepsilon}^{1} +
u^{2}\int\limits_{\varepsilon}^{1}\log (t) e^{-u^{2}t}dt
\\[5mm] & \displaystyle =
\log (\varepsilon)\cdot\left(1-e^{-u^{2}\varepsilon}\right) +
u^{2}\int\limits_{\varepsilon}^{1}\log (t) e^{-u^{2}t}dt.
\end{array}
$$
For fixed $u$, we have that $1-e^{-u^{2}\varepsilon} =
O(\varepsilon)$ as $\varepsilon$ approaches zero, so then
$$
\int\limits_{0}^{1}\left(e^{-u^{2}t}-1\right)\frac{dt}{t} =
u^{2}\int\limits_{0}^{1}\log (t) e^{-u^{2}t}dt.
$$
If we let $v = u^{2}t$, then $dv = u^{2}dt$, hence
$$
\int\limits_{0}^{1}\left(e^{-u^{2}t}-1\right)\frac{dt}{t}=\int\limits_{0}^{u^{2}}\log
(v/u^{2}) e^{-v}dv = \int\limits_{0}^{u^{2}}\log (v) e^{-v}dv -
\log (u^{2})\int\limits_{0}^{u^{2}} e^{-v}dv.
$$
By definition
$$
\Gamma(s) = \int\limits_{0}^{\infty}e^{-v}v^{s}\frac{dv}{v},
$$
which implies that
$$
\Gamma'(1) = \int\limits_{0}^{\infty}e^{-v}\log (v) dv =
\int\limits_{0}^{u^{2}}e^{-v}\log (v) dv +
\int\limits_{u^{2}}^{\infty}e^{-v}\log (v) dv,
$$
thus yielding the relation
\begin{equation}
\begin{array}{ll} \displaystyle
\int\limits_{0}^{1}\left(e^{-u^{2}t}-1\right)\frac{dt}{t}&
\displaystyle =\Gamma'(1) -\int\limits_{u^{2}}^{\infty}e^{-v}\log
(v) dv- \log (u^{2}) \int\limits_{0}^{u^{2}} e^{-v}dv
\label{intbound1}
\\[7mm] &\displaystyle =
\Gamma'(1) - \log(u^{2}) + \log(u^{2})e^{-u^{2}} -
\int\limits_{u^{2}}^{\infty}\log (v) e^{-v}dv.
\end{array}
\end{equation}
For $x > 1$, we have the trivial inequalities $0 \leq \log x \leq
2e^{x/2}$, so then
\begin{equation}
\label{est1} 0 \leq \int\limits_{u^{2}}^{\infty}\log (v) e^{-v}dv
\leq 2\int\limits_{u^{2}}^{\infty} e^{-v/2}dv = 4e^{-u^{2}/2} =
o(1) \,\,\,\,\,\textrm{as $u \rightarrow \infty$.}
\end{equation}
Trivially, $\log(u^{2})e^{-u^{2}} = o(1)$ as $u \rightarrow
\infty$.  By substituting this estimate and (\ref{est1}) into
(\ref{intbound1}), the proof of the lemma is complete. \hfill
$\Box$

\begin{nn}\label{5.7}
\textbf{Corollary.} \emph{With notation as above, we have}
$$
\begin{array}{ll}\displaystyle
\int\limits_{0}^{1}&\displaystyle \left(
\theta_{N(u)}(u^{2}t) -
V(N)\left(e^{-2u^{2}t}I_{0}(2u^{2}t)\right)^{d} - 1+
e^{-u^{2}t}\right)\frac{dt}{t}\\[7mm]& \displaystyle \rightarrow
\int\limits_{0}^{1}\left(\Theta_A(t)-V(A)(4\pi
t)^{-d/2}\right)\frac{dt}{t} + \Gamma'(1) - \log (u^{2})  +
o(1)\,\,\,\,\,\,\,\textrm{as $u \rightarrow \infty$.}
\end{array}
$$
\end{nn}

\proof Simply combine Proposition \ref{5.5} and Lemma \ref{5.6}.
 \hfill $\Box$

\begin{nn}\label{5.8}
\textbf{Theorem.} \emph{With notation as above, we have}
$$
\HN_{N(u)}(0) = \log u^{2} -\zeta_{A}'(0)  + o(1)
\,\,\,\,\,\,\,\textrm{as $u \rightarrow \infty$.}
$$
\emph{Equivalently, we have}
$$
\log \left(\prod\limits_{\Lambda_{j}\neq 0}\Lambda_{j}\right)
= V(N(u))\mathcal{I}_d(0) + \log u^{2} - \zeta_{A}'(0) +
o(1)\,\,\,\,\,\,\,\textrm{as $u \rightarrow \infty$}
$$
\emph{where}
$$
\mathcal{I}_d(0) =
- \int\limits_{0}^{\infty}
\left(e^{-2dt}I_{0}(2t)^{d}-e^{-t}\right)\frac{dt}{t}.
$$
\end{nn}

\proof To begin, combine Proposition \ref{5.4} and Corollary
\ref{5.7} to get the asymptotic expansion of
$V(N(u))\HN_{N(u)}(0)$  out to $o(1)$.  Then combine the terms in
the expansion with the expression in (\ref{zetaderivzero}) to
complete the proof of the first assertion.  One then substitutes
the first assertion into the identity stated in Theorem \ref{3.6}
to prove the second relation.
 \hfill $\Box$

\begin{nn}\label{5.9}
\textbf{Remark.} Using the notation of spectral determinants, one
can re-write the main result in Theorem \ref{5.8} as stating the
asymptotic relation
$$
\mbox{$\log\det^{\ast}$}{\mathbf \Delta}_{DT,N(u)} =V(N(u))  \mathcal{I}_d(0) + \log
u^{2} + \mbox{$\log\det^{\ast}$}{\mathbf \Delta}_{RT,A} +
o(1)\,\,\,\,\,\,\,\textrm{as $u \rightarrow \infty$.}
$$
In words, the asymptotic expansion of the log-determinant of the
Laplacian on the discrete torus $DT_{N(u)}$ has a lead term
$$
\begin{array}{ll}
\displaystyle V(N(u)) \mathcal{I}_d(0) + \log u^{2} & \displaystyle =
-V(N(u))\int\limits_{0}^{\infty}\left(e^{-2dt}(I_{0}(2t))^{d}-e^{-t}\right)\frac{dt}{t}
+ \log u^{2}
\\[7mm] & \displaystyle =
V(N(u))\left(\log 2d -
\int\limits_{0}^{\infty}\left(e^{-dt}(I_{0}(t)^{d}-1\right)\frac{dt}{t}\right)
+ \log u^{2}
\end{array}
$$
with a constant term equal to the log-determinant of the
Laplacian on the real torus $A\mathbf{Z}^{d}\backslash \mathbf{R}^{d}$.
\end{nn}

\begin{nn}\label{5.10}
\textbf{Remark.} As previously discussed, there are two zeta
functions being considered:  The spectral zeta function $\zeta_{A}$ associated to the real torus
$A\mathbf{Z}^{d}\backslash \mathbf{R}^{d}$ and the Epstein zeta
function $\ezeta(\cdot, Q)$ associated to the positive definite $d
\times d$ matrix $Q$.  Classical and elementary mathematical
considerations show that $\ezeta(\cdot, Q)= \zeta_{A}$
where $Q$ is the form associated to $A^{\ast}$, the dual lattice
of $A$.  This remark needs to be kept in mind when reading the
main result in Theorem \ref{5.8} as restated in Remark \ref{5.9}.
\end{nn}

\section{Example: the cases $d=1$ and $d=2$}\label{6}

As stated in the introduction, the lead term in the asymptotic of
$\log\det^{\ast}{\mathbf \Delta}_{DT,N(u)}$ for $d=2$ was first
computed by Kasteleyn \cite{Kasteleyn}. Duplantier and David
\cite{DD} expressed the next order term in terms of the Dedekind
$\eta$ function.  In this section we rederive these results by the
methods introduced in this paper and go further by indicating how
to obtain the complete asymptotic expansion of
$\log\det^{\ast}{\mathbf \Delta}_{DT,N(u)}$ when $d=2.$  We begin
by making explicit the case $d=1$.

\begin{nn}\label{6.0} \textbf{The case $d=1$.}
Known evaluations of the Riemann zeta function
$\zeta_{\mathbf{Q}}$ at $s=0$ yield the relations
$\zeta_{\mathbf{Q}}(0)=-1/2$ and $\zeta_{\mathbf{Q}}'(0) =
-(1/2)\log (2\pi)$.  Let $N(u)=n$ with $n/u \rightarrow \alpha >
0$.  With all this, we can write
\begin{equation}\label{d_equal_one}
\log n^{2} = \log u^{2} +  4 \log
(2\pi/\alpha)\zeta_{\mathbf{Q}}(0) - 4 \zeta_{\mathbf{Q}}'(0)
+o(1) \,\,\,\,\,\,\,\textrm{as $u \rightarrow \infty$.}
\end{equation}
We claim that (\ref{d_equal_one}) is equivalent to the Main
Theorem when $d=1$.  First, it is obvious that the number of
spanning trees in $n\mathbf{Z}\backslash \mathbf{Z}$ is $n$ and,
hence, the determinant of the combinatorial Laplacian is $n^{2}$,
so the left-hand-side of (\ref{maintheorem}) is $\log n^{2}$.
Regarding the right-hand-side, it is possible to directly verify
through the Mellin transform that
$$
\int\limits_{0}^{\infty}e^{-2t}(I_{0}(2t)-1)\frac{dt}{t} = \log 2,
$$
which implies, in our notation, that $\Id_{1}(0)=0$.  Hence, the
lead term in the asymptotic expansion (\ref{maintheorem}) when
$d=1$ is zero.  The continuous Laplacian on $n \mathbf{Z}
\backslash \mathbf{R}$ has the set of eigenvalues given by
$\{(2\pi)^{2}(n/\alpha)^{2}\}$ for $n \in \mathbf{Z}$, so the
spectral zeta function is
$$
\zeta_{\alpha}(s) = 2 (2\pi/\alpha)^{-2s}\zeta_{\mathbf{Q}}(2s),
$$
so
$$
\zeta_{\alpha}'(0)=-\mbox{$\log\det^{\ast}$}{\mathbf
\Delta}_{RT,\alpha} = -4\log(2\pi/\alpha)\zeta_{\mathbf{Q}}(0) + 4
\zeta'_{\mathbf{Q}}(0),
$$
which agrees with (\ref{d_equal_one}) and confirms, by direct
computation, the Main Theorem when $d=1$.
\end{nn}


\begin{nn}\label{6.2}
\textbf{The lead term when $d=2$.} We will evaluate the
expression
$$
\Id_d(0) = \left(  \log (2d) -
\int\limits_{0}^{\infty}\left(e^{-dt}(I_{0}(t)^{d}-1\right)\frac{dt}{t}\right)
$$
when $d=2$. Recall from section (\ref{2.8}) the Mellin inversion
formula, namely equation (\ref{eq:5a}).  Let us set $x=y=0$ and
move the contour of integration to the line $\textrm{Re}(\sigma) =
- 1/4$. The integrand has a pole at $z=0$, which results in the
expression
\begin{equation}
  \label{eq:7}
  \widetilde{f_0f_0}(s)=\frac{1}{2\pi i}\int\limits_{\textrm{Re}(z) = (-1/4)}
\frac{\Gamma(s-z)\Gamma(1/2+z-s)\Gamma(z)\Gamma(1/2-z)} {\pi
\Gamma(1+z-s)\Gamma(1-z)} dz +
\frac{\Gamma(s)\Gamma(1/2-s)\Gamma(1/2)}{\pi\Gamma(1-s)\Gamma(1)},
\end{equation}
so then, using the definition of $\Gamma(s)$ as the Mellin
transform of $e^{-t}$, we have
\begin{equation}
  \label{eq:8}\begin{array}{ll}
\displaystyle \widetilde{f_0f_0}(s)-\Gamma(s)  = &\displaystyle
\int_0^\infty e^{-t}[I_0(t/2)^2 -1]t^s\frac{dt}{t}
\\ [5mm] & \displaystyle
=  \frac{1}{2\pi i}\int\limits_{\textrm{Re}(z) = (-1/4)}
\frac{\Gamma(s-z)\Gamma(1/2+z-s)\Gamma(z)\Gamma(1/2-z)}
{\pi\Gamma(1+z-s)\Gamma(1-z)} dz \\[7mm] &\displaystyle +
\frac{\Gamma(s)\Gamma(1/2-s)}{\sqrt\pi\Gamma(1-s)}-\Gamma(s).
\end{array}
\end{equation}
By Taylor's theorem, we can write
$$
\lim\limits_{s \rightarrow
0}\left(\frac{\Gamma(s)\Gamma(1/2-s)}{\sqrt\pi\Gamma(1-s)}-\Gamma(s)\right)
= \frac{-\Gamma'(1/2)+\sqrt{\pi}\Gamma'(1)}{\sqrt{\pi}},
$$
where we have used that $\Gamma(1/2)=\sqrt{\pi}$ and
$\Gamma(1)=1$.  One now uses the relation
$$
\Gamma(z)\Gamma(z+1/2) = 2^{1-2z}\sqrt{\pi}\Gamma(2z)
$$
and, by computing the derivative at $z=1/2$, one can show that
$$
\frac{-\Gamma'(1/2)+\sqrt{\pi}\Gamma'(1)}{\sqrt{\pi}} = \log 4.
$$
With these preliminary results, we have shown, by taking $s
\rightarrow 0$ in (\ref{eq:8}), that
\begin{equation}
\begin{array}{ll}
  \label{eq:9}\displaystyle
  \int_0^\infty e^{-t}[I_0(t/2)^2 -1]\frac{dt}{t}&\displaystyle =
 \frac{1}{2\pi i}\int\limits_{\textrm{Re}(z) = (-1/4)}
\frac{\Gamma(-z)\Gamma(1/2+z)\Gamma(z)\Gamma(1/2-z)}
{\pi\Gamma(1+z)\Gamma(1-z)} dz+\log 4\\[7mm]
&\displaystyle =\log 4-\frac{1}{2\pi i}\int\limits_{\textrm{Re}(z)
= (-1/4)} \frac{1} {z^2\cos \pi z} dz,
\end{array}
\end{equation}
where the last equality of integrals comes from employing the
identities
$$
\Gamma(x)\Gamma(1-x)=\pi/\sin\pi x \text{\ and\ }
x\Gamma(x)=\Gamma(x+1).
$$

We now evaluate the last integral in (\ref{eq:9}) by moving the
line of integration toward $\textrm{Re}(z) = - \infty$. The
integrand has poles $z=-n-1/2, $ for $n=0,1,2...$ with residues
$(-1)^n (n+1/2)^{-2}/\pi.$  Thus
\begin{equation}
  \label{eq:14}
   \int_0^\infty e^{-t}[I_0(t/2)^2 -1]t^s\frac{dt}{t}=
\log 4-\frac{1}{\pi}\sum_{n=0}^\infty \frac{(-1)^n}{(n+1/2)^2}=
\log 4-\frac {4G}{\pi}
\end{equation}
\end{nn}
where $G$ is the Catalan constant
$$
G = \sum_{n=0}^\infty \frac{(-1)^n}{(2n+1)^2}.
$$
With all this, we have shown that
\begin{equation}\label{firstorderterm}
\Id_2(0) = \frac {4G}{\pi}.
\end{equation}
The relation in (\ref{firstorderterm}), which is the main
asymptotic term from Theorem \ref{5.8} in the case $d=2$, agrees
with the computations from \cite{DD}, specifically equation (3.18)
on page 349; see also equation (A.21) on page 427 of \cite{DD}.

\begin{nn}\label{6.3}
\textbf{Secondary terms when $d=2$.} From Theorem \ref{5.8}, we
have that
\begin{equation}\label{secondterm}
\HN_{N}(0) = \log u^{2} + \mbox{$\log\det^{\ast}$}{\mathbf
\Delta}_{RT,d,A} + o(1) \,\,\,\,\,\,\,\textrm{as
$u \rightarrow 0$}
\end{equation}
for a degenerating sequence of discrete tori.  As before, consider
a sequence of integer vectors $N(u)$ for $u \in \mathbf Z^{+}$ and
assume
$$
\frac{1}{u} N(u) = \frac{1}{u} (n_{1}(u), n_{2}(u)) \rightarrow
(\alpha_{1},\alpha_{2}) \,\,\,\,\,\,\,\textrm{as $u \rightarrow
\infty$.}
$$
Proposition \ref{4.7} evaluates the heat
kernel on $(\alpha_{1}{\mathbf Z} \times \alpha_{2}{\mathbf
Z})\backslash {\mathbf R}^{2}$ where the Laplacian is
$$
\Delta_{\mathbf R^{2}} =-\left(\frac{\partial^{2}}{\partial
x_{1}^{2}} + \frac{\partial^{2}}{\partial x_{2}^{2}}\right)
$$
where $x_{1}$ and $x_{2}$ are the usual global coordinates on
$\mathbf R^{2}$.  The action of $\Delta_{\mathbf R^{2}}$ on the
space of smooth functions on $(\alpha_{1}{\mathbf Z} \times
\alpha_{2}{\mathbf Z})\backslash {\mathbf R}^{2}$ has the set of
eigenvalues given by $\{(2\pi)^{2}((n/\alpha_{1})^{2} +
(m/\alpha_{2})^{2})\}$ for $n, m \in \mathbf Z$, so then its
associated spectral zeta function is
$$
\zeta_{A}(s) = (2\pi)^{-2s} \sum\limits_{(n,m)\neq (0,0)}
\frac{1}{((n/\alpha_{1})^{2} + (m/\alpha_{2})^{2})^{s}}.
$$
Let us write
$$
\zeta_{A}(s) = (\alpha_{1}\alpha_{2})^{s} E(z,s)
$$
where
$$
E(z,s)= (2\pi)^{-2s}\sum\limits_{(n,m)\neq (0,0)}
\frac{\textrm{Im}(z)^{s}}{\vert n z+
m\vert^{2s}}\,\,\,\,\,\textrm{with} \,\,\,\,\, z =
i(\alpha_{2}/\alpha_{1}) = iy.
$$
The function $E(z,s)$ admits a mermorphic continuation to all $s
\in \mathbf C$ with expansion
$$
E(z,s) = -1 - s \log(y\vert \eta(z)\vert^{4}) + O(s)
\,\,\,\,\,\,\,\textrm{as $s \rightarrow 0$}
$$
where $z = i(\alpha_{2}/\alpha_{1}) = iy$ and $\eta(z)$ is the
classical Dedekind eta function
$$
\eta(z) = e^{2\pi i z}\prod\limits_{n=0}^{\infty}\left(1 - e^{2\pi
i n z}\right) \,\,\,\,\,\,\,\, \textrm{for any $z \in \mathbf C$
with $\textrm{Im}(z) > 0$}.
$$
With all this, we then have that
$$\mbox{$\log\det^{\ast}$}
{\mathbf \Delta}_{RT,A} =
-\frac{\partial}{\partial s}\zeta_{A}\Big|_{s=0} =
\log (\alpha_{1}\alpha_{2}) +  \log(y\vert \eta(z)\vert^{4}).
$$
Therefore,
\begin{equation}\label{secondorderterm}
\HN_{N}(0) = \log u^{2} + \log
(\alpha_{1}\alpha_{2}) +  \log(y\vert \eta(z)\vert^{4}) = \log
(n_{1}(u)n_{2}(u)) +  \log(y\vert \eta(z)\vert^{4}).
\end{equation}
The expression (\ref{secondorderterm}) agrees with the
second-order term given in equation (A.21) on page 427 of
\cite{DD}, where those authors used the notation $M=n_{1}(u)$ and
$N = n_{2}(u)$.
\end{nn}

\begin{nn}\label{6.4}
\textbf{Error terms when $d=2$.} The asymptotic expansion in
Theorem \ref{5.8} contains an error term of $o(1)$.  We present
here a technique which improves the error term, ultimately
obtaining an expansion with error term of arbitrarily small
polynomial order as stated in the introduction.  Additionally, the
computations we give here give a second proof of the Main Theorem
in the case $d=2$.

Recall the Mellin inversion formula from section (\ref{2.8}),
namely
\begin{equation}\label{formula1}
\int\limits_{0}^{\infty}I_{x}(t/2)I_{y}(t/2)e^{-t}t^{s}\frac{dt}{t}
= \frac{1}{2\pi i}\int\limits_{\textrm{Re}(z) =
\sigma}\frac{\Gamma(s-z+x)\Gamma(\frac{1}{2}+z-s)\Gamma(z+y)\Gamma(\frac{1}{2}-z)}
{\pi \Gamma(x+1+z-s)\Gamma(y+1-z)}dz.
\end{equation}
For fixed $a$ and $b$, Stirling's formula yields
\begin{equation}\label{stirling}
\frac{\Gamma(z+a)}{\Gamma(z+b)} = z^{a-b} + O(z^{a-b-1})
\end{equation}
as $z$ tends to infinity anywhere in the half-plane
$\textrm{Re}(z) > \delta > 0$.  Substituting (\ref{stirling}) into
(\ref{formula1}), we have that the lead term in the asymptotic
expansion of (\ref{formula1}) is
\begin{equation}\label{formula2}
\begin{array}{ll} \displaystyle
& \displaystyle \frac{1}{2\pi i}\int\limits_{\textrm{Re}(z)
=\sigma} x^{s-z-1-z+s}y^{z-1+z} \cdot
\frac{1}{\pi}\Gamma(z+\frac{1}{2}-s)\Gamma(\frac{1}{2}-z)dz \\[7mm] & \,\,\,\,\,\,\,\,\,
\displaystyle = \frac{x^{2s}}{\pi \cdot xy}\cdot \frac{1}{2\pi
i}\int\limits_{\textrm{Re}(z) =\sigma}
\left(\frac{y}{x}\right)^{2z}\Gamma(z+\frac{1}{2}-s)\Gamma(\frac{1}{2}-z)dz.
\end{array}
\end{equation}

Let us analyze (\ref{formula2}) by moving the contour of
integration toward $\textrm{Re}(z) = \infty$. The integrand has
poles at $z=1/2+n$ with residues equal to $(-1)^{n}/n!$, so then
(\ref{formula2}) is equal to
\begin{equation}\label{formula3}
S = \frac{x^{2s}}{\pi \cdot xy}\cdot \sum\limits_{n=0}^{\infty}
\left(\frac{y}{x}\right)^{2n+1}\frac{(-1)^{n}}{n!}\Gamma(n-s+1).
\end{equation}
Re-write (\ref{formula3}) as
$$
S = \frac{y^{2s-2}}{\pi}\cdot \sum\limits_{n=0}^{\infty}
\left(\frac{y}{x}\right)^{2n-2s+2}\frac{(-1)^{n}}{n!}\Gamma(n-s+1)
$$
and use the elementary formula
\begin{equation}\label{Gammaint}
\Gamma(w)a^{-w}=\int\limits_{0}^{\infty}e^{-at}t^{w}\frac{dt}{t}
\end{equation} with $w=n-s+1$ and $a=x^{2}/y^{2}$ to get that
$$
S = \frac{y^{2s-2}}{\pi}\cdot\sum\limits_{n=0}^{\infty}
\int\limits_{0}^{\infty}e^{-(x^{2}/y^{2})t} t^{n-s+1}\cdot
\frac{(-1)^{n}}{n!}\frac{dt}{t}.
$$
Observing that
$$
e^{-t}=\sum\limits_{n=0}^{\infty}t^{n}\cdot \frac{(-1)^{n}}{n!}
$$
we arrive at the relation
\begin{equation}\label{formula4}
S = \frac{y^{2s-2}}{\pi}\cdot
\int\limits_{0}^{\infty}e^{-(x^{2}/y^{2})t} e^{-t}
t^{-s+1}\frac{dt}{t}.
\end{equation}
The integral can be evaluated using (\ref{Gammaint}) with $w=1-s$
and $a=(x^{2}/y^{2}+1)$, yielding
\begin{equation}\label{formula5}
S = \frac{y^{2s-2}}{\pi}\cdot
\Gamma(1-s)\left(x^{2}/y^{2}+1\right)^{s-1} =\frac{1}{\pi}
\Gamma(1-s)\left(x^{2}+y^{2}\right)^{s-1}.
\end{equation}
If we set $s=-w$ with $\textrm{Re}(w)>0$, we then have that
\begin{equation}\label{leadasymp}
\int\limits_{0}^{\infty}I_{x}(t/2)I_{y}(t/2)e^{-t}t^{-w}\frac{dt}{t}
= \frac{1}{\pi} \cdot
\frac{\Gamma(1+w)}{\left(x^{2}+y^{2}\right)^{1+w}} + \textrm{lower
order terms}.
\end{equation}
Finally, we set $x=nj$ and $y=mk$ with fixed $j$ and $k$ and
positive integers $n$ and $m$ and sum over all $n$ and $m$, thus
proving that
\begin{equation}\label{sumformula}
\sum\limits_{n,m}\int\limits_{0}^{\infty}I_{nj}(t/2)I_{mk}(t/2)e^{-t}t^{-w}\frac{dt}{t}
=\frac{\Gamma(1+w)}{\pi}
\sum\limits_{n,m}\frac{1}{\left((nj)^{2}+(mk)^{2}\right)^{1+w}} +
\textrm{lower order terms}.
\end{equation}
The series in (\ref{sumformula}) can be related to the Eisenstein
series $E_{A}(z,s)$ with $A=(x,y)$ and $z=im/n$.  One then can use the classical
Kronecker's limit formula to evaluate the asymptotic behavior of
(\ref{sumformula}) as $w$ approaches zero.  In all, the above
computations provide another proof of the Main Theorem when $d=2$.

Stirling's formula can be used to compute further terms in the
asymptotic expansion in (\ref{stirling}). Specifically, for any
positive integers
 $N$ and $M$, we write \begin{equation}\label{higherstirling}
\frac{\Gamma(x-z)\Gamma(z+y)}{\Gamma(x+1+z)\Gamma(y+1-z)} -
\frac{1}{xy}\left(\frac{y}{x}\right)^{2z} =
\sum\limits_{j,k}P_{j,k}(z)x^{-2z-1-j}y^{2z-1-k} +
O\left(x^{-2z-1-N}y^{2z-1-M}\right) \end{equation} as $x, y
\rightarrow \infty$, where $P_{j,k}$ denotes a polynomial in $z$,
and the sum is over all integers $j$ and $k$ such that $0 \leq j
\leq N$ and $0 \leq k \leq M$ provided $(j,k) \neq (0,0)$.  By
comparing with (\ref{formula1}), we see that the lower order terms
in (\ref{sumformula}) can be determined by studying
$$
S_{j,k}(x,y)
= \frac{1}{\pi \cdot x^{j+1}y^{k+1}}\cdot \frac{1}{2\pi i}\int\limits_{\textrm{Re}(z)
=\sigma}P_{j,k}(z) \left(\frac{y}{x}\right)^{2z}\Gamma(z+\frac{1}{2})\Gamma(\frac{1}{2}-z)dz. $$
Moving the contour of integration,
we get that
$$ S_{j,k}(x,y) = \frac{1}{\pi \cdot x^{j+1}y^{k+1}}\cdot
\sum\limits_{n=0}^{\infty}(-1)^{n}P_{j,k}(n+1/2) \left(\frac{y}{x}\right)^{2n+1}. $$ Since
$P_{j,k}(z)$ is a polynomial in $z$, we have that $$ P_{j,k}(n+1/2) \left(\frac{y}{x}\right)^{2n+1}
= P_{j,k}\left(\frac{y}{2}\frac{d}{dy}\right)\left(\frac{y}{x}\right)^{2n+1}, $$
so then \begin{equation}
\begin{array}{ll}\label{generalterm} \displaystyle S_{j,k}(x,y) & \displaystyle
= \frac{1}{\pi \cdot x^{j+1}y^{k+1}}\cdot P_{j,k}\left(\frac{y}{2}\frac{d}{dy}\right)
\sum\limits_{n=0}^{\infty} (-1)^{n} \left(\frac{y}{x}\right)^{2n+1}\\[7mm]& \displaystyle = \frac{1}{\pi \cdot x^{j+1}y^{k+1}}\cdot
P_{j,k}\left(\frac{y}{2}\frac{d}{dy}\right) \left(\frac{xy}{x^{2}+y^{2}}\right). \end{array} \end{equation} Let us define
\begin{equation}\label{coefficient1} B_{j,k}(u) = \sum\limits_{(n_{1},n_{2}) \neq (0,0)}S_{j,k}(n_{1}u,n_{2}u).
\end{equation}
Elementary bounds show that $$ B_{j,k}(u) =
O\left(u^{-j-k-2}\right) \,\,\,\,\,\,\,\textrm{when $u \rightarrow
\infty$, provided $(j,k) \neq (0,0)$} $$ since the series in
(\ref{coefficient1}) is convergent for $(j,k) \neq (0,0)$.
Finally, if we set \begin{equation} \label{seriescoefficients}
F_{n}(u) = \sum\limits_{j+k=n}B_{j,k}(u), \end{equation} we have
that the error term $o(1)$ in Theorem \ref{5.8} can be improved to
$$ \sum\limits_{n=1}^{K}F_{n}(u) + O\left(u^{-K-3}\right)
\,\,\,\,\,\textrm{as $u\rightarrow \infty$ for any integer $K >
0$.} $$ To summarize, the coefficients (\ref{seriescoefficients})
can be explicitly determined from the evaluation
(\ref{generalterm}) which utilizes the precise evaluation of
Stirling's formula as stated in (\ref{higherstirling}). The
coefficients (\ref{seriescoefficients}) can be made explicit by
evaluation the polynomials $P_{j,k}(z)$ as defined in
(\ref{higherstirling}).  We will not pursue this analysis further
here, but instead leave the details for consideration elsewhere.
\end{nn}

\section{Additional considerations}\label{7}

In this section we present a number of computations and remarks
discussing further aspects of the preceeding analysis.

\begin{nn}\label{6.8a} \textbf{A precise evaluation of the lead term.}
Let
$$ a_{d} = \log (2d) - \int\limits_{0}^{\infty}e^{-dt}(I_{0}(t)^{d}-1)\frac{dt}{t}
= \log 2 -
\int\limits_{0}^{\infty}\left(e^{-dt}I_{0}(t)^{d}-e^{-t}\right)\frac{dt}{t},
$$
which is obtained by replacing $2t$ by $t$ in the integral
definition of $\Id_{d}(0)$.  With this,
$$ a_{d+1} - a_{d} =
\int\limits_{0}^{\infty}e^{-dt}(I_{0}(t)^{d}(1-e^{-t}I_{0}(t))\frac{dt}{t}.
$$ Since $I_{0}(t) \leq e^{t}$, we have that $a_{d+1} \geq a_{d}$.
Let us write $$ e^{-t}I_{0}(t) = 1 +
\sum\limits_{k=1}^{\infty}b_{k}t^{k}, $$ which implies that
$$ a_{d+1} - a_{d} = \sum\limits_{k=1}^{\infty}b_{n}\int\limits_{0}^{\infty}e^{-dt}(I_{0}(t))^{d}t^{k}dt. $$
From \cite{KN},
we have the identity $$ F(u) =
\int\limits_{0}^{\infty}e^{-ut}I_{0}(t)dt = \frac{1}{u\sqrt{u+2}}.
$$ Therefore, $$
\int\limits_{0}^{\infty}e^{-dt}(I_{0}(t))^{d}t^{k}dt =
(-\partial_{u})^{k}F^{\ast d}(d), $$ where $F^{\ast d}$ denotes
the $d$-fold convolution of $F$ with itself, evaluated at $d$.
Therefore, we have $$ a_{d+1} = a_{d} +
\sum\limits_{k=1}^{\infty}b_{k}\cdot (-\partial_{u})^{k}F^{\ast
d}(d). $$ In section \ref{6.2} we have shown that $a_{2} =
4G/\pi$, where $G$ is the Catalan constant. In addition, one can
write
$$ e^{-t}I_{0}(t) = \frac{1}{\pi}\int\limits_{0}^{\pi}e^{-t\sin^{2}(\theta/2)}d\theta, $$
so the coefficients $b_{k}$ are explicitly computable.
With all this, we have presented another series expansion for
$a_{d}$ which could be used to numerically evaluate $a_{d}$.
\end{nn}

\begin{nn}\label{6.7}
\textbf{Numerical estimation of the lead term.} This subsection
should be compared with the discussion in Felker-Lyons
\cite{FelkerLyons}.
By focusing on the
lead term in Theorem \ref{5.8}, we have that
\begin{equation}\label{riemannsum}
\frac{1}{V(N(u))} \sum\limits_{\Lambda_{j}(u) \neq 0} \log
\Lambda_{j}(u) = \log 2d -
\int\limits_{0}^{\infty}e^{-dt}(I_{0}(t)^{d}-1)\frac{dt}{t} +
o(1)\,\,\,\,\,\textrm{as $u \rightarrow \infty$.}
\end{equation}
As stated in the introduction, the explicit form of the
eigenvalues are such that we recognize the left-hand-side of
(\ref{riemannsum}) as a Riemann sum, so then (\ref{riemannsum})
implies the identity
\begin{equation}\label{leadtermcomparison}
\int\limits_{{\mathbf Z}^{d}\backslash {\mathbf
R}^{d}}\log\left(2d-2\cos(2\pi x_{1})- \cdots - 2\cos(2\pi
x_{d})\right)dx_{1}\cdots dx_{d} =\log 2d -
\int\limits_{0}^{\infty}e^{-dt}(I_{0}(t)^{d}-1)\frac{dt}{t}.
\end{equation}  Recall that
$$
I_{0}(t) = 1+\sum\limits_{n=1}^{\infty}a_{n}t^{2n}
\,\,\,\,\,\textrm{where $a_{n} = \frac{1}{2^{2n}(n!)^{2}}$.}
$$
Let us write
$$
(I_{0}(t))^{d} = 1 + \sum\limits_{n=1}^{\infty}a_{n,d}t^{2n}.
$$
The coefficients $a_{n,d}$ easily can be written in terms of
$a_{n}$ and multinomial coefficients, and
(\ref{leadtermcomparison}) can be evaluated as
$$
\int\limits_{{\mathbf Z}^{d}\backslash {\mathbf
R}^{d}}\log\left(2d-2\cos(2\pi x_{1})- \cdots - 2\cos(2\pi
x_{d})\right)dx_{1}\cdots dx_{d} =\log 2d -
\sum\limits_{n=1}^{\infty}\frac{a_{n,d}\Gamma(2n)}{n^{2d}}.
$$
The numerical evaluation of the left-hand-side of
(\ref{leadtermcomparison}) is difficult for large $d$, both
because of the slow growth of the integrand and because the domain
of integration involves the $d$-fold product of the unit interval.
By comparison, note that the right-hand-side of
(\ref{leadtermcomparison}) allows for rapid numerical evaluation
which improves as $d$ gets larger.  Indeed, from Lemma \ref{4.1},
we have
\begin{equation}\label{seriestail}
0 \leq
\int\limits_{t_{0}}^{\infty}e^{-dt}(I_{0}(t)^{d}-1)\frac{dt}{t}
\leq C^{d}\int\limits_{t_{0}}^{\infty}t^{-d/2}\frac{dt}{t} +
\frac{1}{t_{0}}\int\limits_{t_{0}}^{\infty}e^{-dt}dt =
\frac{2C^{d}t_{0}^{-d/2+1} + e^{-dt_{0}}}{dt_{0}}.
\end{equation}
This bound obviously gets closer to zero as $t_{0}$ grows for any
$d$. In addition, the bound decays exponentially in $d$ provided
$Ct_{0}^{-1/2} < 1$ which, from the estimate that $C <
1/\sqrt{2}$, implies that we need $t_{0} > 1/2$.  By the Mean
Value Theorem, we have, for any $z$ and $a$, the estimate
$$
(x+a)^{d} - x^{d} \leq ad(x+a)^{d-1}.
$$
For any $N$, to be chosen later, we take
$$
x = \sum\limits_{n=0}^{N-1}a_{n}t^{2n}
\,\,\,\,\,\,\,\textrm{and}\,\,\,\,\,\,\, a = I_{0}(t) -x <
I_{0}(t)
$$
to arrive at the bound
$$
\int\limits_{0}^{t_{0}}e^{-dt}(I_{0}(t)^{d}-1)\frac{dt}{t} \leq
\int\limits_{0}^{t_{0}}e^{-dt}\left(\left(\sum\limits_{n=0}^{N-1}a_{n}t^{2n}\right)^{d}-1\right)\frac{dt}{t}
+
d\int\limits_{0}^{t_{0}}e^{t}(e^{-t}I_{0}(t))^{d-1}\left(\sum\limits_{n=N}^{\infty}a_{n}t^{2n}\right)\frac{dt}{t}.
$$
Choose and fixed $t_{0} \in (1/2,1)$.  For any $\varepsilon > 0$,
there is an $N$ such that
$$
\sum\limits_{n=N}^{\infty}a_{n}t^{2n} \leq \varepsilon t^{2N}
\,\,\,\,\,\,\,\textrm{for $t \in [0,t_{0}]$.}
$$
Elementary arguments prove the existence of a constant $A>0$ such
that
$$
e^{-t}I_{0}(t) \leq 1-At \,\,\,\,\,\,\,\textrm{for $t \in
[0,t_{0}]$.}
$$
Therefore,
$$
\int\limits_{0}^{t_{0}}e^{t}(e^{-t}I_{0}(t))^{d-1}\left(\sum\limits_{n=N}^{\infty}a_{n}t^{2n}\right)\frac{dt}{t}
\leq \varepsilon d
e^{t_{0}}\int\limits_{0}^{t_{0}}t^{2N-1}(1-At)^{d-1}dt \leq
\varepsilon e^{t_{0}}t_{0}^{d}/A,
$$
so then
\begin{equation}\label{boundinn}
0 < \int\limits_{0}^{t_{0}}e^{-dt}(I_{0}(t)^{d}-1)\frac{dt}{t} -
\int\limits_{0}^{t_{0}}e^{-dt}\left(\left(\sum\limits_{n=0}^{N-1}a_{n}t^{2n}\right)^{d}-1\right)\frac{dt}{t}
\leq \varepsilon e^{t_{0}}t_{0}^{2N-1}/A.
\end{equation}
Observe that the upper bound in (\ref{boundinn}) is independent of
$d$. When combining (\ref{seriestail}) with (\ref{boundinn}), we
conclude that one can estimate
\begin{equation}\label{termestimate}
\int\limits_{0}^{\infty}e^{-dt}(I_{0}(t)^{d}-1)\frac{dt}{t}\,\,\,\,\,\,\,\textrm{by}\,\,\,\,\,\,\,
\int\limits_{0}^{t_{0}}e^{-dt}\left(\left(\sum\limits_{n=0}^{N-1}a_{n}t^{2n}\right)^{d}-1\right)\frac{dt}{t}
\end{equation}
with error which is independent of $d$.

It remains to study the asymptotic behavior in $d$ of the second
integral in (\ref{termestimate}).  For this, one needs to expand
the $d$-fold product of the polynomial integrand and carry out the
integration.  The integral of each term can then be expressed in
terms of the incomplete Gamma function, which in turn can be
estimated by the Gamma function itself.  One then would employ
elementary but somewhat involved expressions for multinomial
coefficients, ultimately obtaining the following result.  For any
integer $k \geq 1$, there is a polynomial $P_{k}(x)$ of degree $k$
with $P_{k}(0) = 0$ such that
$$
\int\limits_{0}^{t_{0}}e^{-dt}\left(\left(\sum\limits_{n=0}^{N-1}a_{n}t^{2n}\right)^{d}-1\right)\frac{dt}{t}
= P_{k}(1/d) + O(d^{-k-1}) \,\,\,\,\,\,\,\textrm{as $d \rightarrow
\infty$.}
$$
We leave the details of these computations to the interested
reader.
\end{nn}

\begin{nn}\label{6.5}
\textbf{The spectral zeta function.} In very general
circumstances, the spectral zeta function is defined as the Mellin
transform of the theta function formed with the non-zero
eigenvalues. Specifically, for $w \in \mathbf C$ with
$\textrm{Re}(w) > d/2$, we have that
\begin{equation}\label{spectralzeta}
\zeta_{N}(w)={\mathbf M}\theta_{N}(w) =
\frac{1}{\Gamma(w)}\int\limits_{0}^{\infty}\left(\theta_{N}(t)-1\right)t^{w}\frac{dt}{t}.
\end{equation}
As in section 5, let us consider a degenerating sequence of
discrete tori.  The elementary change of variables $t \mapsto
u^{2}t$ in (\ref{spectralzeta}) yields the expression
\begin{equation}\label{spectralzeta2}
\zeta_{N}(w) =
\frac{u^{2w}}{\Gamma(w)}\int\limits_{0}^{\infty}
\left(\theta_{N}(u^{2}t)-1\right)t^{w}\frac{dt}{t}.
\end{equation}
Proposition \ref{5.2} establishes the pointwise convergence of the
integrand in (\ref{spectralzeta2}), and Lemma \ref{5.3} proves a
uniform upper bound.  Combining these results, as in the proof of
Proposition \ref{5.4}, together with (\ref{spectralzeta2}), we
conclude that
\begin{equation}\label{spectralconvergence1}
\lim\limits_{u \rightarrow \infty}\left(u^{-2w}\zeta_{N(u)}(w)\right)
=\zeta_{A}(w)
\end{equation}
for any $w\in \mathbf C$ with $\textrm{Re}(w) > d/2$ and where
\begin{equation}
\label{realspectralzeta} \zeta_{A}(w) =
\frac{1}{\Gamma(w)}\int\limits_{0}^{\infty}
\left(\Theta_A(t)-1\right)t^{w}\frac{dt}{t}
\end{equation}
is the spectral zeta function on the real torus $A\mathbf Z^{d}
\backslash \mathbf Z^{d}$.  The difficulty is determining the
correct generalization of (\ref{spectralconvergence1}) for $w \in
\mathbf C$ with $\textrm{Re}(w)
> 0$.  The spectral zeta function (\ref{spectralzeta2}) can be
meromorphically continued to $w\in \mathbf C$ with $\textrm{Re}(w)
> 0$ by writing
\begin{equation}
\label{spectralzeta3}
\begin{array}{ll} \displaystyle u^{-2w}\zeta_{N}(w) &\displaystyle =
\frac{1}{\Gamma(w)}\int\limits_{1}^{\infty}\left(\theta_{N}(u^{2}t)-1\right)
t^{w}\frac{dt}{t}
\\[7mm] & \displaystyle +
\frac{1}{\Gamma(w)}\int\limits_{0}^{1}\left(\theta_{N(u)}(u^{2}t)
- V(N)e^{-2du^{2}t}(I_{0}(2u^{2}t))^{d}\right)t^{w}\frac{dt}{t}
\\[7mm] & \displaystyle
- \frac{1}{\Gamma(w)}
\int\limits_{0}^{1}\left(V(N)e^{-2du^{2}t}(I_{0}(2u^{2}t))^{d} - 1
\right)t^{w}\frac{dt}{t}.
\end{array}
\end{equation}
Similarly, the meromorphic continuation of
(\ref{realspectralzeta}) is obtained by the expression
\begin{equation}
\label{realspectralzeta2} \begin{array}{ll} \displaystyle
\zeta_{A}(w) &\displaystyle =
\frac{1}{\Gamma(w)}\int\limits_{1}^{\infty}\left(\Theta_A(t)-1\right)t^{w}\frac{dt}{t}
\\[7mm] & \displaystyle +\frac{1}{\Gamma(w)}
\int\limits_{0}^{1}\left(\Theta_A(t)-V(A)(4\pi
t)^{-d/2}\right)t^{w}\frac{dt}{t} \\[7mm] & \displaystyle -
\frac{1}{\Gamma(w)}\int\limits_{0}^{1}\left(V(A)(4\pi
t)^{-d/2}-1\right)t^{w}\frac{dt}{t}.
\end{array}
\end{equation}
The convergence of the first integral in (\ref{spectralzeta3}) to
the first integral in (\ref{realspectralzeta2}) follows from
Proposition \ref{5.2} and Lemma \ref{5.3}, as in the proof of
(\ref{spectralconvergence1}).  The convergence of the second
integral in (\ref{spectralzeta3}) to the second integral in
(\ref{realspectralzeta2}) follows from Proposition \ref{5.5},
using The Lebesgue Dominated Convergence Theorem.  Let us write
the third integral in (\ref{spectralzeta3}) as
\begin{equation}\label{thirddiscrete}
\int\limits_{0}^{1}\left(V(N)e^{-2du^{2}t}I_{0}(2u^{2}t)^{d} - 1
\right)t^{w}\frac{dt}{t}   = \frac{V(N)}{u^{2w}}
\int\limits_{0}^{u^{2}}\left(e^{-2t}I_{0}(2t)\right)^{d}t^{w}\frac{dt}{t}
- \frac{1}{w}
\end{equation}
and the third integral in (\ref{realspectralzeta}) as
\begin{equation}\label{thirdreal}
\int\limits_{0}^{1}\left(V(A)(4\pi
t)^{-d/2}-1\right)t^{w}\frac{dt}{t} =
\frac{V(A)}{(4\pi)^{d/2}(w-d/2)} - \frac{1}{w}.
\end{equation}
Combining all of these results, we have shown that
\begin{equation}
\label{spectralconvergence2} \lim\limits_{u \rightarrow \infty}
u^{-2w}\left(\zeta_{N}(w) -
\frac{V(N)}{\Gamma(w)}\int\limits_{0}^{u^{2}}\left(e^{-2t}(I_{0}(2t))
\right)^{d}t^{w}\frac{dt}{t}\right)
= \zeta_{A}(w)
-\frac{V(A)}{(4\pi)^{d/2}(w-d/2)\Gamma(w)}.
\end{equation}
Using Proposition \ref{4.7}, one can easily show that
(\ref{spectralconvergence2}) is equivalent to
(\ref{spectralconvergence1}) when $\textrm{Re}(w) > d/2$, noting
that (\ref{spectralconvergence2}) holds for all $w \in \mathbf C$
with $\textrm{Re}(w) > 0$ provided $\textrm{Re}(w) \neq d/2$.  By
continuity, (\ref{spectralconvergence2}) extends to all $w \in
\mathbf C$ with $\textrm{Re}(w) > 0$ provided one uses the
interpretation
\begin{equation}\label{highermoment}
\lim\limits_{w \rightarrow d/2}\left( \zeta_{A}(w)
-\frac{V(A)}{(4\pi)^{d/2}(w-d/2)\Gamma(w)}\right) =
\textrm{CT}_{w=d/2}\zeta_{A}(w),
\end{equation}
the constant term in the Laurent expansion at $w=d/2$.

\vskip .10in The above computations in the special case $d=2$ with
$w=1$ and $w=2$ were obtained in section 3 of \cite{DD}.  In the
case when $w=1$ and $d=2$, so then $w=d/2$, the results in
\cite{DD} were expressed in terms of the modular forms.  However,
in light of the general Kronecker limit formula for Epstein zeta
functions and its functional equation, as stated in section
\ref{2.9}, it is evident that for general $d \neq 1$, the limiting
value obtained in (\ref{highermoment}) can be expressed in terms of a
 modular form (\ref{kronecker}).
\end{nn}

\begin{nn}\label{6.6}
\textbf{The Epstein-Hurwitz zeta function.} Further analysis in
\cite{DD} involves the investigation of the finite product
$\prod(s^{2}+\Lambda_{j})$ for degenerating families of
two-dimensional discrete tori for general $s \in \mathbf C$.  If
$s=0$, then Theorem \ref{5.8} determines the asymptotic behavior
of the product of non-zero eigenvalues.  For $s \neq 0$, we one
includes the zero eigenvalues, which trivially introduces the
multiplicative factor of $s^{2}$.  In this setting, we recall that
Theorem \ref{3.6} establishes a relation for
$\prod(s^{2}+\Lambda_{j})$ in terms of integrals involving
$I$-Bessel functions.  Using the substitution $t \mapsto u^{2}t$,
we arrive at the expression
\begin{equation}\label{detexpression}
\log \prod((s/u)^{2}+\Lambda_{j}) = \sum\limits_{\Lambda_{j} \neq
0}\log \left((s/u)^{2}+\Lambda_{j}\right) = V(N/u)\Id_d(s/u) +
\HN_{N}(s/u).
\end{equation}
where
$$
\Id_d(s/u)=
-\int\limits_{0}^{\infty}\left(e^{-2du^{2}t}e^{-s^{2}t}(I_{0}(2u^{2}t))^{d}-e^{-u^{2}t}\right)\frac{dt}{t}
$$
and
$$
\HN_{N}(s/u)  = - \int\limits_{0}^{\infty}
e^{-s^{2}t}\left(\theta_{N}(u^{2}t) -V(N)
e^{-2du^{2}t}(I_{0}(2u^{2}t))^{d}\right) \frac{dt}{t} - \log(s^{2}).
$$
Trivially, one has that
$$
\prod((s/u)^{2}+\Lambda_{j}) = u^{-2V(N)}\prod(s^{2}+u^{2}\cdot
\Lambda_{j}).
$$
We can now employ bounds from section 5 in order to determine the
asymptotic behavior of (\ref{detexpression}) as $u \rightarrow
\infty$.  Specifically, one uses Proposition \ref{4.7},
Proposition \ref{5.2}, Lemma \ref{5.3}, and Proposition \ref{5.5}
to show that for any $s \in \mathbf C$ with $\textrm{Re}(s^{2})>
0$ we have the asymptotic formula
$$
\prod((s/u)^{2}+\Lambda_{j}) = V(A/u)\Id_d(s/u) + \HN_{A}(s) +o(1)
\,\,\,\,\,\,\,\textrm{as $u \rightarrow \infty$}
$$
where
$$
\HN_{A}(s)  = - \int\limits_{0}^{\infty}
e^{-s^{2}t}\left(\Theta_A(t) -V(A) (4\pi t)^{-d/2}\right)
\frac{dt}{t} - \log(s^{2}).
$$
The function $\HN_{A}(s)$ is related to the regularized
harmonic series, as defined and studied in \cite{JLang0}
associated to the set $\{\Lambda_{A,j}+s^{2}\}$ where
$\{\Lambda_{A,j}\}$ is the set of eigenvalues of the Laplacian on
the real torus $A\mathbf Z^{d}\backslash {\mathbf R}^{d}$. The
general results from \cite{JLang0} establish that
$\HN_{A}(s)$ can be expressed as special values of the
Epstein-Hurwitz zeta function formed with the set
$\{\Lambda_{A,j}+s^{2}\}$.  We refer the reader to \cite{JLang0}
for further details and identities.
\end{nn}

\begin{nn}\label{6.8}
\textbf{General discrete tori.} Let $B$ be a positive definite $d
\times d$ integer matrix and consider the discrete tori $DT_{B}
= B{\mathbf Z}^{d} \backslash {\mathbf Z}^{d}$.  The results from
\cite{KN} easily extend to compute the spectrum of the Laplacian
on $DT_{B}$ in terms of the dual lattice $B^{\ast}$. The
existence and uniqueness of the associated heat kernel on
$DT_{B}$ follows from general results (see, for example,
\cite{Dodziuk} and \cite{DM}), thus allowing one to extend the
results of section 3 above.  Assume there is a one-parameter
family $DT_{B(u)}$ of discrete tori parameterized by $u \in
\mathbf Z$ such that $B(u)/u \rightarrow M$ as $u\rightarrow
\infty$ where $M$ is a  positive definite $d \times d$ matrix. The
results of the present article apply when considering spectral
invariants on $DT_{B(u)}$.  For example, Theorem \ref{5.8} will
extend to obtain the asymptotic behavior of the determinant of the
Laplacian on $DT_{B(u)}$ with second-order term equal to the
zeta-regularized spectral determinant on the real torus $M{\mathbf
Z}^{d}\backslash {\mathbf R}^{d}$.  Additional results of the
present paper, specifically the contents of section \ref{6.5} and
section \ref{6.6}, carry through using the proofs given and only a
slight change
in notation.
\end{nn}

\begin{nn}\label{6.10}
\textbf{Height functions.} The asymptotic expansion from Theorem
\ref{5.8}, and its extension as outlined in section \ref{6.8}, can
be interpreted as saying that $\mbox{$\log\det^{\ast}$}{\mathbf
\Delta}_{DT,\ast}$ is a height function on the space of discrete
tori of fixed dimension since it tends to $+\infty$ through
degeneration.  In the case of real tori, one needs to introduce a
minus sign and study $-\mbox{$\log\det^{\ast}$}{\mathbf
\Delta}_{RT,\ast}$ in order to have a height function which tends
to $+\infty$ through degeneration;  see, for example, \cite{Chiu}.
The problem of finding real tori with minimum height remains a
question of interest, see, for example, \cite{SaSt}. One point of
future investigation is to see to what extent Theorem \ref{5.8}
allows for a connection between the problems considered in
\cite{SaSt} and the study of the height function
$\mbox{$\log\det^{\ast}$}{\mathbf \Delta}_{DT,\ast}$ on discrete
tori, which has the advantage of being an invariant which is
defined as a finite product and not through meromorphic
continuation.

\end{nn}


\vspace{5mm}


\noindent
Gautam Chinta \\
Department of Mathematics \\
The City College of New York \\
Convent Avenue at 138th Street \\
New York, NY 10031
U.S.A. \\
e-mail: chinta@sci.ccny.cuny.edu

\vspace{5mm} \noindent
Jay Jorgenson \\
Department of Mathematics \\
The City College of New York \\
Convent Avenue at 138th Street \\
New York, NY 10031
U.S.A. \\
e-mail: jjorgenson@mindspring.com

\vspace{5mm}\noindent
Anders Karlsson \\
Mathematics Department \\
Royal Institute of Technology \\
100 44 Stockholm, Sweden \\
e-mail: akarl@kth.se

\end{document}